\documentclass[a4paper,12pt, oneside, reqno]{amsart}
\usepackage{mathrsfs}
\usepackage{amsfonts}
\usepackage{amsmath,amstext,amssymb,amsthm}
\usepackage[polish, english]{babel}
\usepackage{amsxtra}
\usepackage{color}
\usepackage{dsfont}
\usepackage[margin=2cm, centering]{geometry}
\usepackage[bbgreekl]{mathbbol}
\usepackage{soul}
\usepackage{graphics}

\usepackage[colorlinks,citecolor=blue,urlcolor=blue,bookmarks=true]{hyperref}

\DeclareSymbolFontAlphabet{\mathbb}{AMSb}
\DeclareSymbolFontAlphabet{\mathbbl}{bbold}

\usepackage{accents}

\theoremstyle{plain}
\newtheorem{theorem}{Theorem}[section] 
\newtheorem{lemma}[theorem]{Lemma}
\newtheorem{proposition}[theorem]{Proposition}
\newtheorem{corollary}[theorem]{Corollary}
\newtheorem{fact}[theorem]{Fact}
 \theoremstyle{definition}
\newtheorem{defn}{Definition}

\newtheorem{example}{Example} 
\newtheorem*{remark*}{Remark} 
\newtheorem{remark}[theorem]{Remark}

\newcommand{\R}{\mathbb{R}}
\newcommand{\Rd}{{\R^{d}}}
\newcommand{\NN}{\mathbb{N}}
\newcommand{\ind}{\mathds{1}}
\renewcommand{\leq}{\leqslant}
\renewcommand{\geq}{\geqslant}

\def\ind{{\bf 1}}
\def\qed{{\hfill $\Box$ \bigskip}}
\def\LL{{\mathcal L}}
\def\RR{{\mathbb R}}

\def\N{{\mathbb N}}
\def\pf{\noindent{\bf Proof.} }

\def\({\left(} 
\def\){\right)} 
\def\[{\left[}
\def\]{\right]} 
\def\<{\langle} 
\def\>{\rangle}

\def \Qzero{\rm{(Q0)}}
\def \Qa{\rm{(Q1)}}
\def \Qb{\rm{(Q2)}}

\def \Ab{\rm{(A2)}}

\newcommand{\lah}{\alpha_h}
\newcommand{\uah}{\beta_h}

\newcommand{\err}[2]{\rho_{#1}^{#2}} 
\newcommand{\param}{\sigma}  
\newcommand{\lmCJ}{c_{\!\scriptscriptstyle J}}  

\newcommand{\aF}{\mathcal{F}} 

\newcommand{\drf}{b} 

\newcommand{\rr}{\Upsilon} 

\newcommand{\ka}{\kappa_3}  
\newcommand{\kb}{\kappa_4}  
\newcommand{\TCh}{\Theta} 
\newcommand{\ct}{\vartheta} 

\definecolor{ks}{rgb}{0.7,0.1,0.2}

\definecolor{zm}{RGB}{128,128,0}

\allowdisplaybreaks

\title{Fundamental solution for super-critical non-symmetric L{\'e}vy-type operators}

\thanks{The research was partially supported by
 the 
National Science Centre (Poland)
grant 2016/23/B/ST1/01665.}

\author[K. Szczypkowski]{Karol Szczypkowski}
\address{
	Karol Szczypkowski\\
	Wydzia{\lll} Matematyki,
	Politechnika Wroc{\lll}awska\\
	Wyb. Wyspia\'{n}skiego 27\\
	50-370 Wroc{\lll}aw\\
	Poland
}
\email{karol.szczypkowski@pwr.edu.pl}

\date{}

\begin{document}

\begin{abstract}
We prove the existence and give estimates of the fundamental solution (the heat kernel) for the equation $\partial_t =\LL^{\kappa}$ for non-symmetric non-local operators
$$
\LL^{\kappa}f(x):= \int_{\R^d}( f(x+z)-f(x)- \ind_{|z|<1} \left<z,\nabla f(x)\right>)\kappa(x,z)J(z)\, dz\,,
$$
under broad assumptions on $\kappa$ and $J$.
Of special interest is the case when the order of the operator $\LL^{\kappa}$ is smaller than or equal to 1. Our approach rests on imposing
suitable cancellation conditions on the internal drift coefficient
$$
\int_{r\leq |z|<1} z \kappa(x,z)J(z)dz\,,\qquad 0<r\leq 1\,,
$$
which allows us to handle the non-symmetry of
$z\mapsto \kappa(x,z)J(z)$.
The results are new even for the $1$-stable L{\'e}vy measure $J(z)=|z|^{-d-1}$. 
\end{abstract}

\maketitle



\section{Introduction}\label{sec:intr}

In recent years, there has been a lot of interest in constructing
semigroups for L{\'e}vy-type operators
\cite{MR3353627},
\cite{MR3652202},
\cite{MR3500272},
\cite{MR3817130},
\cite{GS-2018},
\cite{MR2163294}, 
\cite{MR2456894},
\cite{FK-2017},
\cite{BKS-2017},
\cite{KR-2017},
\cite{MR3294616},
\cite{MR3544166},
\cite{PJ}, 
\cite{CZ-new},
\cite{CZ-survey},
\cite{MR3965398}.
Such operators arise naturally due to the Courr{\`e}ge-Waldenfels theorem \cite[Theorem~4.5.21]{MR1873235}, \cite[Theorem~2.21]{MR3156646}.
In general, they are not symmetric, so the $L^2$-theory or Dirichlet forms 
and the corresponding $L^2$-semigroups of operators
does not apply in this context.
We shall discuss operators of the form
\begin{align}
\LL^{\kappa}f(x)&:= \int_{\Rd}( f(x+z)-f(x)- \ind_{|z|<1} \left<z,\nabla f(x)\right>)\kappa(x,z)J(z)\, dz \,, \label{e:intro-operator-a1-crit1}
\end{align}
and allow for non-symmetric measures  $\kappa(x,z)J(z) dz$. 
The paper is a continuation of the 
research conducted in \cite{GS-2018}, where
operators of the form \eqref{e:intro-operator-a1-crit1} were considered under stronger conditions.
We also improve and extend the results of \cite{PJ}, and part of those in~\cite{CZ-new} and
\cite{MR3652202}.

We now introduce our setting, notation and motivations.
Let $d\in\N$ and
$\nu:[0,\infty)\to[0,\infty]$ be a non-increasing  function ($\nu \not\equiv 0$) with
$$
\int_{\Rd}  (1\land |x|^2) \nu(|x|)dx<\infty\,.
$$
We consider
$J: \Rd  \to [0, \infty]$ 
 such that for some  $\lmCJ \in [1,\infty)$ and 
 all $x\in \Rd$,
\begin{equation}\label{e:psi1}
\lmCJ^{-1} \nu(|x|)\leq J(x) \leq \lmCJ \,\nu(|x|)\,.
\end{equation}
Furthermore, suppose that 
$\kappa(x,z)$ is a Borel 
function on $\R^d\times \Rd$ such that
\begin{equation}\label{e:intro-kappa}
0<\kappa_0\leq \kappa(x,z)\leq \kappa_1<\infty\, , 
\end{equation}
for some numbers $\kappa_0$, $\kappa_1$,
and there is $\beta\in (0,1)$ 
and a number $\kappa_2\geq 0$
with
\begin{equation}\label{e:intro-kappa-holder}
|\kappa(x,z)-\kappa(y,z)|\leq \kappa_2|x-y|^{\beta}\, .
\end{equation}
The following concentration functions  play a prominent role in the paper,
$$
h(r):= \int_{\Rd} \left(1\land \frac{|x|^2}{r^2}\right) \nu(|x|)dx\,,\qquad \quad
K(r):=r^{-2} \int_{|x|<r}|x|^2 \nu(|x|)dx\,,\qquad r>0\,.
$$
We say that
\emph{the weak  scaling condition} at the origin holds if there are $\lah\in (0,2]$ and $C_h \in [1,\infty)$ such that 
\begin{equation}\label{eq:intro:wlsc}
 h(r)\leq C_h\,\lambda^{\lah}\,h(\lambda r)\, ,\qquad  r, \lambda \in (0,1]\,.
\end{equation}
In a  similar fashion, we consider the existence of $\uah\in (0,2]$ and $c_h\in (0,1]$ such that
\begin{equation}\label{eq:intro:wusc}
 h(r)\geq c_h\,\lambda^{\uah}\,h(\lambda r)\, ,\qquad r, \lambda \in (0,1]\,.\\
\end{equation}
We propose two more conditions on $\kappa(x,z)$ and $J(z)$.
Suppose there are  numbers $\ka, \kb\geq 0$ such that
\begin{align}\label{e:intro-kappa-crit}
\sup_{x\in\Rd}\left| \int_{r\leq |z|<1} z\, \kappa(x,z) J(z)dz \right| &\leq \ka  rh(r)\,,
 \qquad r\in (0,1],\\
\left| \int_{r\leq |z|<1} z\, \big[ \kappa(x,z)- \kappa(y,z)\big] J(z)dz \right| &\leq \kb |x-y|^{\beta} rh(r)\,, \qquad r\in (0,1]. \label{e:intro-kappa-crit-H}
\end{align}

\vspace{\baselineskip}
\noindent
We are ready to specify our framework.
The {\it dimension} $d$ and {\it the profile function $\nu$} are fixed.
We will use 
alternatively
 two sets of assumptions.
\begin{enumerate}
\item[] 
\begin{enumerate}
\item[$\Qa$:] \quad \eqref{e:psi1}--\eqref{eq:intro:wlsc}
hold, $\lah=1$; \eqref{e:intro-kappa-crit} and \eqref{e:intro-kappa-crit-H} hold;
\item[$\Qb$:] \quad \eqref{e:psi1}--\eqref{eq:intro:wusc}
hold,  $0<\lah \leq \uah <1$ and $1-\lah<\beta \land \lah$; \eqref{e:intro-kappa-crit} and \eqref{e:intro-kappa-crit-H} hold.
\end{enumerate}
\end{enumerate}

\vspace{\baselineskip}

For the sake of the discussion, we assume
that \eqref{e:psi1}--\eqref{eq:intro:wlsc} hold.
We notify that the internal non-symmetry of the operator \eqref{e:intro-operator-a1-crit1} may result in a non-zero 
{\it internal drift} coefficient 
\begin{align*}
\int_{\Rd} z \left( \ind_{|z|<r}-\ind_{|z|<1} \right) \kappa(x,z)J(z) dz\,.
\end{align*}
Accordingly,
\begin{align}\label{eq:L_split}
\LL^{\kappa}f(x)= 
\int_{\Rd}( f(x+z)-f(x)- \ind_{|z|<r} \left<z,\nabla f(x)\right>)\,\kappa(x,z)J(z) dz \nonumber \\
+\left(\int_{\Rd} z \left( \ind_{|z|<r}-\ind_{|z|<1} \right) \kappa(x,z)J(z) dz\right) \cdot \nabla f(x)\,.
\end{align}
The influence of the {\it internal drift}
term \eqref{eq:L_split}
may be different depending on
the order of the operator $\LL^{\kappa}$
measured by $\lah$ and $\uah$ in \eqref{eq:intro:wlsc} and~\eqref{eq:intro:wusc}.
Below we will usually let $r=h^{-1}(1/t)$.
\begin{fact}\label{fact:1}
If  \eqref{e:psi1}--\eqref{eq:intro:wlsc}
 hold and $\lah>1$, then
\eqref{e:intro-kappa-crit} and \eqref{e:intro-kappa-crit-H} hold.
\end{fact}
\noindent
The fact follows from Lemma~\ref{lem:int_J}.
Thus, we have that if $\lah>1$ (the sub-critical case), then the inequalities
\eqref{e:intro-kappa-crit} and \eqref{e:intro-kappa-crit-H} are automatically satisfied.
It explains a posteriori the success of the analysis of the sub-critical non-symmetric case in \cite{GS-2018}, see also \cite{MR3652202}, \cite{PJ}, \cite{CZ-new}, \cite{MR1744782}.
On the other hand, if  $\lah=1$ (the critical case)  or $\lah< 1$  (the super-critical case), the study of the operator \eqref{e:intro-operator-a1-crit1} is harder and we need the conditions \eqref{e:intro-kappa-crit} and \eqref{e:intro-kappa-crit-H}.
This resolves a question about natural conditions for critical and sub-critical non-symmetric operators that lead to strong results, similar to those in \cite{GS-2018},
and which provide
a thorough analysis of the operator,
the fundamental solution for the corresponding parabolic equation,
 and the associated semigroup. We obtain such results if either of the set of assumptions $\Qa$ or $\Qb$ is satisfied, see Section~\ref{sec:main_res}.
Note that under the symmetry condition, i.e., 
when $J(z)=J(-z)$ and $\kappa(x,z)=\kappa(x,-z)$, $x,z\in\Rd$, 
which is usually assumed in the literature,
the problematic terms involving $\nabla f$ disappear after rewriting the operator $\LL^{\kappa}$ as
\begin{align*}
\frac12 
\int_{\Rd}( f(x+z)+f(x-z)-2f(x))\,\kappa(x,z)J(z)\, dz \,.
\end{align*}
In fact, under the symmetry, we have
\begin{align}\label{eq:0}
\sup_{r\in (0,1]} \sup_{x\in\Rd}\left| \int_{r\leq |z|<1} z \kappa(x,z)J(z)dz \right| =0\,.
\end{align}
The condition \eqref{eq:0} was used in a non-symmetric case 
in \cite{PJ} and \cite{CZ-new} for $\nu(r)=r^{-d-1}$. The conditions \eqref{e:intro-kappa-crit} and
\eqref{e:intro-kappa-crit-H} are much less restrictive.

\begin{example}\label{ex:1}
Let $\nu(r)=r^{-d-1}$.
Then \eqref{eq:intro:wlsc}
and
\eqref{eq:intro:wusc}
hold with $\lah=\uah=1$.
The inequalities \eqref{e:intro-kappa-crit} and
\eqref{e:intro-kappa-crit-H} are read as 
\begin{align*}
 \sup_{x\in\Rd}\left| \int_{r\leq |z|<1} z\, \kappa(x,z) J(z)dz \right| &\leq c \,,
 \qquad r\in (0,1],\\
\left| \int_{r\leq |z|<1} z\, \big[ \kappa(x,z)- \kappa(y,z)\big] J(z)dz \right| &\leq c |x-y|^{\beta}\,, \qquad r\in (0,1]. 
\end{align*}
Hence, if $J(z)$ and $\kappa(x,z)$
are such that \eqref{e:psi1}, \eqref{e:intro-kappa},
\eqref{e:intro-kappa-holder}
and the above inequalities hold, then
the assumptions $\Qa$ are satisfied.
We also note that
$\int_{r \leq|z|<1}|z| \nu(|z|)dz= c \log(1/r)$.
\end{example}

Due to our general setting, we can deal with other interesting operators.

\begin{example}
Let $\nu(r)=r^{-d-1}\log(2+1/r)$.
Then
\eqref{eq:intro:wlsc}
holds with $\lah=1$, but not with any $\lah>1$. 
Furthermore,
\eqref{eq:intro:wusc} holds for  every $\uah>1$,
but not with $\uah=1$.
We also have that $\nu(r)$ is comparable to $r^{-d}h(r)$,
see
\cite[Lemma~5.3 and~5.4]{GS-2018}.
Thus
\eqref{e:intro-kappa-crit} and
\eqref{e:intro-kappa-crit-H}
allow,  respectively, logarithmic growth as $r\to 0$ as follows
\begin{align*}
 \sup_{x\in\Rd}\left| \int_{r\leq |z|<1} z\, \kappa(x,z) J(z)dz \right| &\leq c \log(2+1/r) \,,
 \qquad r\in (0,1],\\
\left| \int_{r\leq |z|<1} z\, \big[ \kappa(x,z)- \kappa(y,z)\big] J(z)dz \right| &\leq c  |x-y|^{\beta} \log(2+1/r)\,, \qquad r\in (0,1]. 
\end{align*}
Thus, if $J(z)$ and $\kappa(x,z)$
are such that \eqref{e:psi1}, \eqref{e:intro-kappa},
\eqref{e:intro-kappa-holder}
and the above inequalities hold, then
the assumptions $\Qa$ are satisfied.
Noteworthy, here $\int_{r \leq|z|<1}|z|\nu(|z|)dz$ is comparable to \mbox{$[\log(2+1/r)]^2$} for small $r$.
\end{example}

\begin{example}
Let $\nu(r)=r^{-d-\alpha}$ and $\alpha \in (1/2,1)$.  Then \eqref{eq:intro:wlsc}
and
\eqref{eq:intro:wusc}
hold with $\lah=\uah=\alpha$.
Note that for $J(z)$ and $\kappa(x,z)$
to satisfy $\Qb$, we need {\it the balance condition} $\alpha+\beta>1$ to hold.
Furthermore, we have $r h(r)=r^{1-\alpha}h(1)$,
while 
$\int_{r \leq|z|<1}|z|\nu(|z|)dz$ is comparable to a positive constant for small $r$.
\end{example}

\noindent
Given a profile function $\nu$, it is not hard to find  $J(z)$ and $\kappa(x,z)$ such that \eqref{e:psi1}--\eqref{e:intro-kappa-holder} hold, and
\begin{align*}
\sup_{x\in\Rd}\left| \int_{r\leq |z|<1} z \kappa(x,z)J(z)dz \right|
\geq 
c \int_{r\leq |z|<1} |z|\nu(|z|)\,dz\,,\qquad r\in (0,1],
\end{align*}
for some $c>0$.
Therefore,
in each of the above examples,
such choice of $J(z)$ and $\kappa(x,z)$
is not admissible, because
the condition
\eqref{e:intro-kappa-crit} fails.
In fact, they can be chosen so that
\eqref{e:intro-kappa-crit-H} fails as well.
Put differently, 
similarly to \eqref{eq:0},
the conditions 
\eqref{e:intro-kappa-crit}
and
\eqref{e:intro-kappa-crit-H}
require certain cancellations to take place.

The success 
of our approach 
based on the usage of
\eqref{e:intro-kappa-crit} and
\eqref{e:intro-kappa-crit-H} suggests that
also in other studies and contexts
where a
counterpart of \eqref{eq:0}
plays a role (see \cite{CZ-Dini})
a  relaxation of assumptions
to proper counterparts of 
\eqref{e:intro-kappa-crit} and
\eqref{e:intro-kappa-crit-H}
should be possible.

As further applications, our results allow solving uniquely  {\it the martingale problem} 
for the operator $(\LL^{\kappa}, C_c^{\infty}(\Rd))$. They also have applications to {\it the Kato class} of the semigroup $(P^{\kappa}_t)_{t\geq 0}$ corresponding to $\LL^{\kappa}$, given in \eqref{e:intro-semigroup}. For details see
\cite[Remarks~1.5 and~1.6]{GS-2018}.

Throughout the paper we make an effort to control 
how constants depend on the initial parameters of our model. The reason for that is twofold.
First of all, 
it is 
necessary in the preliminaries,
like Sections~\ref{sec:analysis_LL} and~\ref{sec:analysis_LL_2}, to be able to execute the construction 
and to find the key properties of
a candidate for the solution in Section~\ref{sec:constr}.
The second reason is more application-oriented: uniform results for families of operators or processes are desired in such areas as
mixing property, multiscale models, homogenization, stationary distribution,
see \cite{MR2779833}, \cite{MR1988467}, \cite{MR4238225}.
The operators we consider resemble those investigated in the study of mean field games \cite{MR4309434}.

The main tool used in this paper is the parametrix method,
proposed by E. Levi \cite{zbMATH02644101} to solve elliptic Cauchy problems.
It was successfully applied
in the theory of partial differential equations 
\cite{zbMATH02629782}, 
\cite{MR1545225},
\cite{MR0003340},
\cite{zbMATH03022319}, with an
overview in the monograph \cite{MR0181836},
as well as in the theory of pseudo-differential operators \cite{MR2093219},
\cite{MR3817130}, \cite{MR3652202}, \cite{FK-2017}, \cite{MR3294616}. 
In particular, operators comparable in a sense to the fractional Laplacian were intensively studied 
by this method
\cite{MR0492880},
\cite{MR616459},
\cite{MR972089},
\cite{MR1744782},
\cite{MR2093219},
also very recently
\cite{MR3500272},
\cite{PJ},
\cite{CZ-new},
\cite{KR-2017}.
More detailed historical comments on the development of the method can be found in \cite[Bibliographical Remarks]{MR0181836} and in the introductions of \cite{MR3652202} and \cite{BKS-2017}.

Basically we follow the scheme  of \cite{GS-2018},
which in turn extends and strengthens
 \cite{MR3817130} and \cite{MR3500272}.
The results in the present paper are of the same type  as in \cite{GS-2018}
with the main progress being the recognition and usage of the conditions \eqref{e:intro-kappa-crit} and \eqref{e:intro-kappa-crit-H}.

Other related papers treat, for instance,
(symmetric)
singular L{\'e}vy measures \cite{BKS-2017}, \cite{KR-2017} or 
 (symmetric) exponential L{\'e}vy measures \cite{KL-2018}.
We also list some papers that
use different techniques to associate a semigroup with an operator
by
symbolic calculus 
\cite{MR0367492}, \cite{MR0499861}, \cite{MR666870},
\cite{MR1659620}, \cite{MR1254818}, \cite{MR1917230},
\cite{MR2163294}, \cite{MR2456894},
Dirichlet forms \cite{MR2778606}, \cite{MR898496}, \cite{MR2492992}, \cite{MR2443765}, \cite{MR2806700}
or perturbation series \cite{MR1310558}, \cite{MR2283957}, \cite{MR2643799}, \cite{MR2876511}, \cite{MR3550165}, \cite{MR3295773}.
For probabilistic methods and applications, we refer the reader to \cite{MR3022725}, 
\cite{MR3544166}, \cite{MR1341116},
 \cite{MR3765882}, \cite{K-2015}, \cite{KR-2017}.

As stated in the abstract, when the present paper was first made public on \href{https://arxiv.org/abs/1807.04257v1}{arXiv:1807.04257v1}, the results were new even for the operators discussed in Example~\ref{ex:1}. 
Those operators are now included in the recent paper \cite{KKS-2020}.

\section{Main results}\label{sec:main_res}

We start by giving an exact meaning to \eqref{e:intro-operator-a1-crit1}.
We apply the operator \eqref{e:intro-operator-a1-crit1}, in a strong or weak sense, only when it is well defined according to 
the following definitions.
Let $f\colon \Rd\to \R$ be a Borel measurable function.

\begin{defn}[\textbf{Strong operator}]
We say that the operator $\LL^{\kappa}f$ is well defined if 
the gradient $\nabla f(x)$ exists and
the corresponding integral 
in \eqref{e:intro-operator-a1-crit1} converges absolutely for every $x\in\Rd$.
\end{defn}

We denote by
$\LL^{\kappa,\varepsilon}f$ 
the expression
\eqref{e:intro-operator-a1-crit1}
with $J(z)$ replaced  
by 
$
J(z)\ind_{|z|>\varepsilon}$, $\varepsilon \in [0,1]$.

\begin{defn}[\textbf{Weak operator}]
We let
\begin{equation*}
\LL^{\kappa,0^+}f(x):=\lim_{\varepsilon \to 0^+}\LL^{\kappa,\varepsilon}f(x)\,,
\end{equation*}
if
the (strong) operators $\LL^{\kappa,\varepsilon}f$ are well defined for $\varepsilon \in (0,1]$,
and the limit 
exists for every~$x\in\Rd$.
\end{defn}

It is clear that the
operator $\LL^{\kappa,0^+}$ is an extension of $\LL^{\kappa,0}= \LL^{\kappa}$, meaning
that if $\LL^{\kappa}f$ is well defined, then so is $\LL^{\kappa,0^+}f$ and $\LL^{\kappa,0^+}f=\LL^{\kappa}f$.
Therefore, it is desired to prove the existence of of solutions to the equation $\partial_t=\LL^{\kappa}$ and the uniqueness of a solution to $\partial_t=\LL^{\kappa,0^+}$.

Here are our main results.
\begin{theorem}\label{t:intro-main}
Assume $\Qa$ or $\Qb$. Let $T>0$.
There is a unique function $p^{\kappa}(t,x,y)$ 
on $(0,T]\times \Rd\times \Rd$ 
such~that
\begin{itemize}
\item[(i)]  For all $t\in(0,T]$, $x,y\in \Rd$, $x\neq y$,
\begin{equation}\label{e:intro-main-1}
\partial_t p^{\kappa}(t,x,y)=\LL_x^{\kappa,0^+}p^{\kappa}(t,x, y)\,.
\end{equation}
\item[(ii)] The function $p^{\kappa}(t,x,y)$ is jointly continuous on $(0,T]\times \Rd\times \Rd$
and
for any
$f\in C_c^{\infty}(\Rd)$,
\begin{equation}\label{e:intro-main-5}
\lim_{t\to 0^+}\sup_{x\in \Rd}\left| \int_{\Rd}p^{\kappa}(t,x,y)f(y)\, dy-f(x)\right|=0\, .
\end{equation}
\noindent
\item[(iii)]
 For every $t_0\in (0,T)$ there are $c>0$ and $f_0\in L^{1}(\Rd)$ 
such that for all $t\in (t_0,T]$, $x,y\in\Rd$,
\begin{equation}\label{e:intro-main-2}
|p^{\kappa}(t,x,y)|\le c f_0(x-y)\,,
\end{equation}
and
\begin{equation}\label{e:intro-main-4}
|\LL_x^{\kappa, \varepsilon}p^{\kappa}(t,x,y)|\leq c \,,\qquad \varepsilon \in (0,1]\,.
\end{equation}
\item[(iv)]  For every $t\in (0,T]$ there is $c>0$ such that for all $x,y\in\Rd$,
\begin{equation}\label{e:intro-main-a1}
|\nabla_x p^{\kappa}(t,x,y)|\leq c\,. 
\end{equation}
\end{itemize}
\end{theorem}

In the next theorem, we collect more qualitative properties of $p^{\kappa}(t,x,y)$.
To this end, 
for $t>0$ and $x\in \R^d$ we define {\it the bound function},
\begin{equation}\label{e:intro-rho-def}
\rr_t(x):=\left( [h^{-1}(1/t)]^{-d}\land \frac{tK(|x|)}{|x|^{d}} \right) .
\end{equation}
It is an integrable function, which may provide sharp estimates for the heat kernel, extending the well known two-sided bounds 
$t^{-d/\alpha}\land t/|x|^{d+\alpha}$
of the fundamental solution to  
$\partial_t = \Delta^{\alpha/2}$, where 
$\Delta^{\alpha/2}:= -(-\Delta^{\alpha/2})$ is
the fractional Laplacian, see
\cite[Thorem~1.1]{GS-2017} as well as more detailed discussion provided in \cite[Section~5]{GS-2017}.
Properties of the bound function can be found also in 
\cite[Section~5]{GS-2018}.

\begin{theorem}\label{t:intro-further-properties}
Assume $\Qa$ or $\Qb$. The following hold true.
\begin{enumerate}
\item[\rm (1)] (Non-negativity) The function $p^{\kappa}(t,x,y)$ is non-negative on $(0,\infty)\times\Rd\times\Rd$.
\item[\rm (2)] (Conservativeness) For all $t>0$, $x\in\Rd$, 
\begin{equation*}
\int_{\Rd}p^{\kappa}(t,x,y) dy =1\, .
\end{equation*}
\item[\rm (3)] (Chapman-Kolmogorov equation) For all $s,t > 0$, $x,y\in \R^d$,
\begin{equation*}
\int_{\R^d}p^{\kappa}(t,x,z)p^{\kappa}(s,z,y)\, dz =p^{\kappa}(t+s,x,y)\, .
\end{equation*}
\item[\rm (4)] (Upper estimate) For every $T>0$ there is $c>0$ such that for all $t\in (0,T]$, $x,y\in \Rd$,
\begin{equation*}
p^{\kappa}(t,x,y) \leq c \rr_t(y-x)\, .
\end{equation*}
\item[\rm (5)] (Fractional derivative) For every $T>0$ there is $c>0$ such that for all
$t\in (0,T]$, $x,y\in\Rd$,
\begin{align*}
|\LL_x^{\kappa } p^{\kappa}(t, x, y)|\leq c t^{-1}\rr_t(y-x)\,.
\end{align*}
\item[\rm (6)] (Gradient)
For every
$T>0$ there is $c>0$ such that for all
$t\in (0,T]$, $x,y\in\Rd$, 
\begin{equation*}
\left|\nabla_x p^{\kappa}(t,x,y)\right|\leq  c\! \left[h^{-1}(1/t)\right]^{-1} \rr_t(y-x)\,. 
\end{equation*}
\item[\rm (7)] (Continuity) The function
$\LL_x^{\kappa} p^{\kappa}(t,x,y)$ is jointly continuous on $(0,\infty)\times \Rd\times\Rd$.
\item[\rm (8)] (Strong operator)
For all $t>0$, $x,y\in\Rd$,
\begin{equation*}
\partial_t p^{\kappa}(t,x,y)= \LL_x^{\kappa}\, p^{\kappa}(t,x,y)\,.
\end{equation*}
\item[\rm (9)] (H\"older continuity) For all $T>0$, $\gamma \in [0,1] \cap[0,\lah)$,
there is  $c>0$ such that for all $t\in (0,T]$ and $x,x',y\in \Rd$,
\begin{equation*}
\left|p^{\kappa}(t,x,y)-p^{\kappa}(t,x',y)\right| \leq c 
 (|x-x'|^{\gamma}\land 1) \left[h^{-1}(1/t)\right]^{-\gamma} \big( \rr_t(y-x)+ \rr_t(y-x') \big).
\end{equation*}
\item[\rm (10)] (H\"older continuity)
For all $T>0$, $\gamma \in [0,\beta)\cap [0,\lah)$,
there is $c>0$ such that for all $t\in (0,T]$ and $x,y,y'\in \Rd$,
\begin{equation*}
\left|p^{\kappa}(t,x,y)-p^{\kappa}(t,x,y')\right| \leq c 
(|y-y'|^{\gamma}\land 1) \left[h^{-1}(1/t)\right]^{-\gamma}  \big( \rr_t(y-x)+ \rr_t(y'-x) \big).
\end{equation*}
\end{enumerate}
The constants in {\rm (4) -- (6)} may be chosen to depend only on $d, \lmCJ, \kappa_0, \kappa_1, \kappa_2, \ka, \kb, \beta,  \lah,  C_h, h, T$.
The same holds for {\rm (9)} and {\rm (10)} but with additional dependence on $\gamma$. 
\end{theorem}

For $t>0$, we define 
\begin{equation}\label{e:intro-semigroup}
P_t^{\kappa}f(x)=\int_{\Rd} p^{\kappa}(t,x,y)f(y)\, dy\, ,\quad x\in \Rd\, ,
\end{equation}
whenever the integral exists in the Lebesgue sense.
We also put $P_0^{\kappa}=\mathrm{Id}$, the identity operator.

\begin{theorem}\label{thm:onC0Lp}
Assume  $\Qa$ or $\Qb$. The following hold true.
\begin{enumerate}
\item[\rm (1)]  $(P^{\kappa}_t)_{t\geq 0}$ is an analytic strongly continuous positive contraction semigroup 
on \mbox{$(C_0(\Rd),\|\cdot\|_{\infty})$.}
\item[\rm (2)]  $(P^{\kappa}_t)_{t\geq 0}$ is an analytic strongly continuous  semigroup on every $(L^p(\Rd),\|\cdot\|_p)$, \mbox{$p\in [1,\infty)$.}
\item[\rm (3)] Let $(\mathcal{A}^{\kappa},D(\mathcal{A}^{\kappa}))$ be the 
generator of $(P_t^{\kappa})_{t\geq 0}$ on $(C_0(\Rd),\|\cdot\|_{\infty})$.\\
 Then
\begin{enumerate}
\item[\rm (a)] $C_0^2(\Rd) \subseteq D(\mathcal{A}^{\kappa})$ and $\mathcal{A}^{\kappa}=\LL^{\kappa}$ on $C_0^2(\Rd)$,
\item[\rm (b)] $(\mathcal{A}^{\kappa},D(\mathcal{A}^{\kappa}))$ is the closure of $(\LL^{\kappa}, C_c^{\infty}(\Rd))$,
\item[\rm (c)] the function $x\mapsto p^{\kappa}(t,x,y)$ belongs to $D(\mathcal{A}^{\kappa})$ for all $t>0$, $y\in\Rd$, and
$$
\mathcal{A}^{\kappa}_x\, p^{\kappa}(t,x,y)= \LL_x^{\kappa}\, p^{\kappa}(t,x,y)=\partial_t p^{\kappa}(t,x,y)\,,\qquad x\in\Rd\,.
$$
\end{enumerate}
\item[\rm{(4)}]  Let $(\mathcal{A}^{\kappa},D(\mathcal{A}^{\kappa}))$ be the 
generator of $(P_t^{\kappa})_{t\geq 0}$ on $(L^p(\Rd),\|\cdot\|_p)$, $p\in [1,\infty)$.\\
 Then
\begin{enumerate}
\item[\rm (a)] $C_c^2(\Rd) \subseteq D(\mathcal{A}^{\kappa})$ and $\mathcal{A}^{\kappa}=\LL^{\kappa}$ on $C_c^2(\Rd)$,
\item[\rm (b)] $(\mathcal{A}^{\kappa},D(\mathcal{A}^{\kappa}))$ is the closure of $(\LL^{\kappa}, C_c^{\infty}(\Rd))$,
\item[\rm (c)] the function $x\mapsto p^{\kappa}(t,x,y)$ belongs to $D(\mathcal{A}^{\kappa})$ for all $t>0$, $y\in\Rd$, and in $L^p(\Rd)$,
$$
\mathcal{A}^{\kappa} \, p^{\kappa}(t,\cdot,y)= \LL^{\kappa}\, p^{\kappa}(t,\cdot,y)=\partial_t p^{\kappa}(t,\cdot,y)\,.
$$
\end{enumerate}
\end{enumerate}
\end{theorem}

Finally, we provide
a lower bound for the heat kernel $p^{\kappa}(t,x,y)$.
It is quite typical that one first proves a lower bound in terms of $h^{-1}$ and $\nu$,
as we do in the first two statements of Theorem~\ref{thm:lower-bound}, cf.
\cite[Remark~5.7 and Section~4]{GS-2017},
\cite[Section~5]{MR4140542}.
In our setting, we have
that $\nu(r) \leq c\, r^{-d}K(r)$, see \cite[Lemma~7.1]{GS-2017}, 
which indicates a possible difference between those lower bounds and the upper bound by $\rr_t(y-x)$.
The converse inequality $\nu(r) \geq c r^{-d}K(r)$ is well understood, see
 \cite[Lemma~7.3]{GS-2017},
and leads to sharp two-sided bounds in the third statement of the theorem.
For abbreviation, we write $\varpi$
to denote the collection of $d, \lmCJ,\kappa_0,\kappa_1,\kappa_2,\ka,\kb,\beta,\lah, C_h, h$.

\begin{theorem} \label{thm:lower-bound}
Assume $\Qa$ or $\Qb$. The following hold true.
\begin{itemize}
\item[(i)] There are $T_0=T_0(\nu,\varpi)>0$ and $c=c(\nu,\varpi)>0$ such that for all $t\in (0,T_0]$, $x,y\in\Rd$,
\begin{equation}\label{e:intro-main-11}
p^{\kappa}(t,x,y)\geq c\left(
[h^{-1}(1/t)]^{-d}\wedge t \nu \left( |x-y|\right)\right).
\end{equation}
\item[(ii)] If additionally $\nu$ is positive, then for every $T>0$ there is $c=c(T,\nu,\varpi)>0$ such that \eqref{e:intro-main-11} holds for $t\in(0,T]$ and $x,y\in\Rd$.\\
\item[(iii)]
If additionally there are $\bar{\beta}\in [0,2)$ and $\bar{c}>0$ such that
$\bar{c} \lambda^{d+\bar{\beta}} \nu (\lambda r) \leq \nu(r)$,
$\lambda \leq 1$, $r>0$, then  for every $T >0$
there is $c=c(T,\nu,\bar{c},\bar{\beta},\varpi)>0$ such that for all $t\in(0,T]$ and $x,y\in\Rd$,
 \begin{equation}\label{e:intro-main-111}
p^{\kappa}(t,x,y)  \geq  c  \rr_t(y-x)\,.
\end{equation}
\end{itemize}
\end{theorem}

\begin{remark}\label{rem:smaller_beta}
If 
\eqref{e:intro-kappa}, \eqref{e:intro-kappa-holder} hold, then $|\kappa(x,z)-\kappa(y,z)|\leq (2\kappa_1 \vee \kappa_2)|x-y|^{\beta_1}$ for every $\beta_1 \in [0,\beta]$.
\end{remark}

\noindent
According to the parametrix method,
the fundamental solution $p^{\kappa}$ is expected to 
be given by 
\begin{align*}
p^{\kappa}(t,x,y)=
p^{\mathfrak{K}_y}(t,x,y)+\int_0^t \int_{\Rd}p^{\mathfrak{K}_z}(t-s,x,z)q(s,z,y)\, dzds\,,
\end{align*}
where $q(t,x,y)$ solves the equation
\begin{align*}
q(t,x,y)=q_0(t,x,y)+\int_0^t \int_{\Rd}q_0(t-s,x,z)q(s,z,y)\, dzds\,,
\end{align*}
and 
$$q_0(t,x,y)=\big(\LL_x^{{\mathfrak K}_x}-\LL_x^{{\mathfrak K}_y}\big) p^{\mathfrak{K}_y}(t,x,y)\,.$$
Here $p^{\mathfrak{K}_w}$ is the heat kernel corresponding to the L{\'e}vy operator $\LL^{\mathfrak{K}_w}$
obtained from the operator $\LL^{\kappa}$ by freezing its coefficients: $\mathfrak{K}_w(z)=\kappa(w,z)$.
In our setting, we draw the initial knowledge on $p^{{\mathfrak K}_w}$ from 
\cite{GS-2017}, which we then
exploit in Sections~\ref{sec:analysis_LL} and~\ref{sec:analysis_LL_2}
to establish further properties.
Already in this preliminary part 
we essentially incorporate \eqref{e:intro-kappa-crit} and
\eqref{e:intro-kappa-crit-H}, which differs 
from \cite{GS-2018}. We also see the effect of the internal drift and the fact that the order of the operator does not have to be strictly larger than one, e.g., Proposition~\ref{thm:delta_crit},
\eqref{e:delta-difference-abs-crit1}, Lemma~\ref{l:some-estimates-3b-crit1}.
In Section~\ref{sec:q} we carry out the construction of  $q$ and so also $p^{\kappa}$.
In view of future developments,
the following remark is notable.
\begin{remark}
We emphasize that the construction of $p^{\kappa}$ is possible, and many preliminary facts hold  
true under a weaker assumption
\begin{enumerate}
\item[]
\begin{enumerate}
\item[$\Qzero$:] \quad \eqref{e:psi1}--\eqref{eq:intro:wlsc}
hold, $\lah\in (0,1]$; \eqref{e:intro-kappa-crit} and \eqref{e:intro-kappa-crit-H} hold.
\end{enumerate}
\end{enumerate}
In particular, see
Lemma~\ref{l:estimates-q0-crit1}, Theorem~\ref{t:definition-of-q-crit1}, Lemma~\ref{lem:phi_cont_xy-crit1},
Proposition~\ref{l:phi-y-abs-cont-crit1}
 and \eqref{e:p-kappa}, \eqref{e:def-phi-y-2}.
\end{remark}

\noindent
The subsequent non-trivial step is to verify that $p^{\kappa}$
is the actual solution.
To this end, in Section~\ref{sec:phi} we need extra
constraints which eventually result in  $\Qa$ and $\Qb$, see for instance
Lemma~\ref{e:L-on-phi-y-crit1}
and the comments preceding Lemma~\ref{lem:phi_pomoc-crit1} and Lemma~\ref{lem:some-est_gen_phi_xy-crit1}.
In Section~\ref{sec:p_kappa} we collect the initial properties of $p^{\kappa}$.
In Section~\ref{sec:Main} we establish a nonlocal maximum principle, analyze the semigroup $(P_t^{\kappa})_{t\geq 0}$, complement the fundamental properties of $p^{\kappa}$, and prove Theorems~\ref{t:intro-main}--\ref{thm:lower-bound}.
Section~\ref{sec:appA} contains auxiliary results.
We give a final comment on the connection with \cite{GS-2018}.
\begin{remark}
The structure of the present paper is similar to that of \cite{GS-2018}  to keep the same train of thought, but also to facilitate the transition between the papers while comparing and identifying the corresponding results. The reason for doing the latter is that the
proofs that are the same as in 
\cite{GS-2018} are reduced to a minimum, we only list which facts are needed, occasionally give general ideas, and refer the reader to \cite{GS-2018} for details.  
We deliberately focus on and explain those aspects that are different from \cite{GS-2018}.
We believe that such a presentation makes the content more comprehensible.
In Lemma~\ref{lem:p-kappa-final-prop-crit1}
we also give a correction of a part of the proof of
\cite[Lemma~4.10]{GS-2018}.
\end{remark}

In what follows, 
the function $\nu$ and
the constants
$d$, $\lmCJ$, $\kappa_0$, $\kappa_1$, $\kappa_2$, $\beta$, 
$\ka$, $\kb$,
$\lah$,  $C_h$, $\uah$, $c_h$
can be regarded as fixed.
Apart from Sections~\ref{sec:Main} and~\ref{sec:appA}, we explicitly formulate 
all assumptions in lemmas, corollaries, propositions, and theorems.
On the other hand,
{\bf in Section~\ref{sec:Main} we assume that either
$\Qa$ or $\Qb$ holds}.

\section{Notation}\label{sec:notation}

For the reader's convenience, we collect inhere the notation  repeatedly used in the paper. 
By $c(d,\ldots)$ we denote a
positive number that depends only on the listed parameters $d,\ldots$. 
By
$\param$ we represent the collection of
$\lmCJ,\kappa_0,\kappa_1,\ka,\lah, C_h, h$.
Throughout the article,
$\omega_d=2\pi^{d/2}/\Gamma(d/2)$ is the surface measure of the unit sphere in $\Rd$.
We use ``$:=$" to denote the definition.
As usual, $a\land b:=\min\{a,b\}$ and $a\vee b := \max\{a,b\}$.

The operator $\LL^{\kappa}$ is given in \eqref{e:intro-operator-a1-crit1}.
The functions 
$h(r)$, $K(r)$, and $\rr_t(x)$
were introduced in Section~\ref{sec:intr} and Section~\ref{sec:main_res}.
Here is a glossary of symbols to be used 
(and explained) below.
For 
$$\mathfrak{K}\colon \Rd \to [0,\infty)\,,$$
we introduce the operator
\begin{align}\label{op:aux}
\LL^{\mathfrak{K}}f(x):=
\int_{\Rd}( f(x+z)-f(x)- \ind_{|z|<1} \left<z,\nabla f(x)\right>)\,\mathfrak{K}(z)J(z)\, dz \,.
\end{align}
The corresponding heat kernel is denoted by
\begin{align}\label{heat_kernel:aux}
p^{\mathfrak{K}}(t,x,y)=p^{\mathfrak{K}}(t,y-x)\,.
\end{align}
We let
\begin{align}\label{e:delta-f-def}
\delta_{1.r}^{\mathfrak{K}} (t,x,y;z)&:=p^{\mathfrak{K}}(t,x+z,y)-p^{\mathfrak{K}}(t,x,y)-\ind_{|z|<r}\left< z,\nabla_x p^{\mathfrak{K}}(t,x,y)\right>,
\end{align}
and
\begin{align*}
\delta^{\mathfrak{K}}(t,x,y;z)&:=\delta_{1.1}^{\mathfrak{K}}(t,x,y;z)\,.
\end{align*}
Thus for $\mathfrak{K}_1$ and $\mathfrak{K}_2$ we have 
\begin{align}\label{eq:L_delta}
\LL_x^{\mathfrak{K}_1} \,p^{\mathfrak{K}_2}(t,x,y)
&=\int_{\Rd}\delta^{\mathfrak{K}_2} (t,x,y;z)\, \mathfrak{K}_1(z)J(z)dz\,,
\end{align}
and
\begin{equation}\label{eq:L_delta_r}
\begin{aligned}
\LL_x^{\mathfrak{K}_1} \,p^{\mathfrak{K}_2}(t,x,y)
&=\int_{\Rd}\delta_{1.r}^{\mathfrak{K}_2} (t,x,y;z)\, \mathfrak{K}_1(z)J(z)dz  \\
&\quad +\left(\int_{\Rd} z \left( \ind_{|z|<r}-\ind_{|z|<1} \right) \mathfrak{K}_1(z) J(z)\, dz\right) \cdot \nabla_x p^{\mathfrak{K}_2}(t,x,y)
\,.  
\end{aligned}
\end{equation}
Starting from Section~\ref{sec:analysis_LL_2},
we shall  use
\begin{align}\label{def:k-frozen}
\mathfrak{K}_w(z):=\kappa(w,z)\,,
\end{align}
which defines $\LL^{\mathfrak{K}_w}f(x)$,  
$p^{\mathfrak{K}_w}(t,x,y)$ and
$\delta_{1.r}^{\mathfrak{K}_w}(t,x,y;z)$.
The main objects in the paper are
\begin{align}\label{e:q0-definition}
q_0(t,x,y):=
\big(\LL_x^{{\mathfrak K}_x}-\LL_x^{{\mathfrak K}_y}\big) p^{\mathfrak{K}_y}(t,x,y)
=
 \int_{\Rd}\delta^{\mathfrak{K}_y}(t,x,y;z)\left(\kappa(x,z)-\kappa(y,z)\right)J(z)dz \,,
\end{align}
\begin{align}\label{e:qn-definition}
q_n(t,x,y):=\int_0^t \int_{\Rd}q_0(t-s,x,z)q_{n-1}(s,z,y)\, dzds\,,
\end{align}
\begin{align}\label{def:q}
q(t,x,y):=\sum_{n=0}^{\infty}q_n(t,x,y)\,,
\end{align}
and
\begin{align}\label{e:p-kappa}
p^{\kappa}(t,x,y):=p^{\mathfrak{K}_y}(t,x,y)+\int_0^t \int_{\Rd}p^{\mathfrak{K}_z}(t-s,x,z)q(s,z,y)\, dzds\,.
\end{align}
The integral part of \eqref{e:p-kappa}
is of special interest and to investigate its properties we introduce
\begin{align}\label{e:phi-y-def}
\phi_y(t,x,s):=\int_{\Rd} p^{\mathfrak{K}_z}(t-s,x,z)q(s,z,y)\, dz\,,
\end{align}
and
\begin{align}\label{e:def-phi-y-2}
\phi_y(t,x):=\int_0^t \phi_y(t,x,s)\, ds =\int_0^t \int_{\Rd}p^{\mathfrak{K}_z}(t-s,x,z)q(s,z,y)\, dzds\, .
\end{align}
Our estimates shall be frequently presented by means of
\begin{align}\label{def:err}
\err{\gamma}{\beta}(t,x):= \left[h^{-1}(1/t)\right]^{\gamma} \left(|x|^{\beta}\land 1\right) t^{-1} \rr_t(x)\,.
\end{align} 
To shorten the notation in Section~\ref{sec:analysis_LL} we shall use the  expressions
\begin{align*}
\aF_{1}&:=\rr_t(y-x-z)\ind_{|z|\geq h^{-1}(1/t)}+ \left[ \left(\frac{|z|}{h^{-1}(1/t)} \right)^2 \land \left(\frac{|z|}{h^{-1}(1/t)} \right) \right]  \rr_t(y-x),\\
\aF_{2}&:=\rr_t(y-x-z)\ind_{|z|\geq h^{-1}(1/t)}+ \left[ \left(\frac{|z|}{h^{-1}(1/t)}\right)\wedge 1\right] \rr_t(y-x).
\end{align*}
Thus, $\aF_1=\aF_1(t,x,y;z)$ and $\aF_2=\aF_2(t,x,y;z)$.
We shall also need the non-increasing function
$$
\TCh(t):= 1+\ln\left(1 \vee \left[h^{-1}(1/t)\right]^{-1}\right),\qquad t>0\,.
$$

We use the following function spaces:
$L^p(\Rd)$ denotes the Lebesgue space with $p\in [1,\infty)$,
$C(D)$ are continuous functions on $D\subseteq \RR^n$, $n\in \NN$.
Furthermore, $C_b(\Rd)$,  $C_0(\Rd)$, $C_c(\Rd)$ are subsets of $C(\Rd)$
of functions that are  bounded,
 vanish at infinity,
and have compact support, respectively.
We write
$f\in C^k(\Rd)$ if
the function and all its derivatives up to (including if finite) order $k\in \NN\cup \{\infty\}$ are elements of $C(\Rd)$;
we similarly understand $C_b^k(\Rd)$,  $C_0^k(\Rd)$, $C_c^k(\Rd)$.
In particular,
$C_c^{\infty}(\Rd)$ are smooth functions with compact support.
The set
$C^{k,\eta}(\Rd)$
consists of functions in $C^{k}(\Rd)$
such that all the derivatives of order $k$ are (uniformly) H{\"o}lder continuous with exponent $0<\eta<1$; we
similarly define
$C_b^{k,\eta}(\Rd)$,  $C_0^{k,\eta}(\Rd)$, $C_c^{k,\eta}(\Rd)$.

\section{Analysis of the  heat kernel of $\LL^{\mathfrak{K}}$}
\label{sec:analysis_LL}

In this section, we consider $\LL^{\mathfrak{K}}$  given by \eqref{op:aux}
with $J(z)$ satisfying \eqref{e:psi1} and \eqref{eq:intro:wlsc},
and
a function~$\mathfrak{K}(z)$ such that
\begin{align}\label{ineq:k-bounded}
0<\kappa_0 \leq \mathfrak{K}(z) \leq \kappa_1\,,
\end{align}
and
\begin{align}\label{ineq:k-int_control}
\left| \int_{r\leq |z|<1} z\, \mathfrak{K}(z) J(z)dz \right|\leq \ka rh(r)\,, \qquad r\in (0,1].
\end{align}
The operator $\LL^{\mathfrak{K}}f$ 
is well defined for
functions $f\in C_c^{\infty}(\Rd)$ and uniquely determines a L{\'e}vy process
and its transition density
$p^{\mathfrak{K}}(t,x,y)$ as represented in \eqref{heat_kernel:aux}. 
Then for all $t>0$, $x,y\in\Rd$,
\begin{equation}\label{eq:p_gen_klas}
\partial_t p^{\mathfrak{K}}(t,x,y)= \LL_x^{\mathfrak{K}}\, p^{\mathfrak{K}}(t,x,y)\,.
\end{equation}
For more information we refer the reader to \cite[Section~6]{GS-2018};
in particular, the condition \cite[(96)]{GS-2018} is satisfied, see \eqref{e:psi1}, 
\cite[(86)]{GS-2018}
 and \eqref{eq:intro:wlsc}.

Clearly, \eqref{ineq:k-bounded}
 corresponds to
\eqref{e:intro-kappa}, while
\eqref{ineq:k-int_control}
corresponds to \eqref{e:intro-kappa-crit}. 
We want to emphasize the role of \eqref{e:intro-kappa-crit} and \eqref{ineq:k-int_control}. In particular, \eqref{ineq:k-int_control} yields the following fundamental upper bound for $p^{\mathfrak{K}}$ and its derivatives.

\begin{proposition}\label{prop:gen_est_crit}
Assume \eqref{e:psi1}, \eqref{eq:intro:wlsc}, \eqref{ineq:k-bounded}, \eqref{ineq:k-int_control}. For every $T>0$ and $\bbbeta\in \mathbb{N}_0^d$ there exists a constant $c=c(d,T,\bbbeta,\param)$
such that for all $t\in (0,T]$, $x,y\in\Rd$,
\begin{align*}
|\partial_x^{\bbbeta} p^{\mathfrak{K}}\left(t,x,y\right)|\leq 
c \left[h^{-1}(1/t) \right]^{-|\bbbeta|} \rr_t(y-x)\,.
\end{align*}
\end{proposition}
\pf
The result follows from \cite[Proposition~5.4 ii)]{GS-2017} with $r_*=1$.
\qed

Here is a lower bound for $p^{\mathfrak{K}}$.

\begin{lemma}\label{prop:gen_est_low-crit1}
Assume \eqref{e:psi1}, \eqref{eq:intro:wlsc}, \eqref{ineq:k-bounded}, \eqref{ineq:k-int_control}.
 For every $T,\theta>0$ there exists a constant $\tilde{c}=\tilde{c}(d,T,\theta,\nu,\param)$
such that for all $t\in (0,T]$ and $|x-y|\leq \theta h^{-1}(1/t)$,
\begin{align*}
p^{\mathfrak{K}}\left(t,x,y\right)\geq \tilde{c} \left[ h^{-1}(1/t)\right]^{-d}\,.
\end{align*}
\end{lemma}
\pf
We use
\cite[Corollary~5.5]{GS-2017}
with $x-y- t\drf_{[h_0^{-1}(1/t)]}$ in place of $x$,
which is allowed since by  \eqref{ineq:k-int_control} we have
$|t\drf_{[h_0^{-1}(1/t)]}|\leq a h_0^{-1}(1/t)$ for $a=a(d,T,\param)$. 
\qed

Proposition~\ref{prop:gen_est_crit}
enables analysis of the increments of the heat kernel.

\begin{lemma}\label{lem:pk-collected}
Assume \eqref{e:psi1}, \eqref{eq:intro:wlsc}, \eqref{ineq:k-bounded}, \eqref{ineq:k-int_control}.
 For every $T>0$ there exists a constant $c=c(d,T,\param)$
such that for all $r>0$, $t\in (0,T]$, $x,x',y,z\in\Rd$ we have
\begin{align}
\left|p^{\mathfrak{K}}(t,x+z,y)-p^{\mathfrak{K}}(t,x,y)\right|&\leq c\, \aF_2(t,x,y;z)\,,\label{ineq:est_diff_1} \\
\left|\nabla_x p^{\mathfrak{K}}(t,x+z,y)-\nabla_x p^{\mathfrak{K}}(t,x,y)\right|&\leq c \left[h^{-1}(1/t)\right]^{-1} \aF_2(t,x,y;z)\,,\label{ineq:est_grad_1}\\
|\delta_{1.r}^{\mathfrak{K}}(t,x,y;z)| &\leq c \big(
\aF_{1}(t,x,y;z)\ind_{|z|<r}+\aF_{2}(t,x,y;z)\ind_{|z|\geq r}\big)\,, \label{ineq:est_delta_1_crit}
\end{align}
and whenever  $|x'-x|<h^{-1}(1/t)$, then
\begin{align}\label{ineq:diff_delta_1_crit}
|\delta_{1.r}^{\mathfrak{K}}(t,x',y;z)-\delta_{1.r}^{\mathfrak{K}}(t,x,y;z)|
\leq c\left(\frac{|x'-x|}{h^{-1}(1/t)}\right)  \big(
\aF_{1}(t,x,y;z)\ind_{|z|<r}+\aF_{2}(t,x,y;z)\ind_{|z|\geq r}\big)\,.
\end{align}
\end{lemma}
\pf The inequalities follow from Proposition~\ref{prop:gen_est_crit} and \cite[Corollary~5.10]{GS-2018}, cf. \cite[Lemma 2.3 -- 2.8]{GS-2018}. The idea of the proof is to
represent the differences as integrals of derivatives
in all cases when the absolute value of the argument increment is smaller than $h^{-1}(1/t)$.
\qed

Due to \cite[Corollary~5.10]{GS-2018}, the inequality \eqref{ineq:est_diff_1} can be  written equivalently as
\begin{align}\label{ineq:est_diff_1*}
\left|p^{\mathfrak{K}}(t,x',y)-p^{\mathfrak{K}}(t,x,y)\right|&\leq c \left(\frac{|x'-x|}{h^{-1}(1/t)} \land 1\right) \big( \rr_t(y-x') + \rr_t(y-x)\big)\,.
\end{align}
The form \eqref{ineq:est_diff_1} is useful for estimating integrals, whereas
\eqref{ineq:est_diff_1*} easily yields what follows.

\begin{lemma}\label{lem:pkw_holder}
Assume \eqref{e:psi1}, \eqref{eq:intro:wlsc}, \eqref{ineq:k-bounded}, \eqref{ineq:k-int_control}.
For every $T>0$ there exists a constant
$c=c(d,T,\param)$
such that for all $t\in(0,T]$, $x,x',y,w \in \Rd$ and $\gamma\in [0,1]$,
\begin{align*}
|p^{\mathfrak{K}}(t,x',y)-p^{\mathfrak{K}}(t,x,y) | 
\leq c (|x-x'|^{\gamma}\land 1) \left[h^{-1}(1/t)\right]^{-\gamma} 
 \big( \rr_t(y-x') + \rr_t(y-x)\big).
\end{align*}
\end{lemma}
In the next result we  estimate
$\LL_x^{\mathfrak{K}_1} p^{\mathfrak{K}_2}(t,x,y)$, crucially using \eqref{ineq:k-int_control}.

\begin{lemma}\label{lem:Lkp_abs}
Assume \eqref{e:psi1}, \eqref{eq:intro:wlsc} and let $\mathfrak{K}_1$, $\mathfrak{K}_2$ satisfy \eqref{ineq:k-bounded}, \eqref{ineq:k-int_control}. For every $T>0$ there exists a constant $c=c(d,T,\param)$ such that
for all $t\in (0,T]$, $x,y\in\Rd$ we have
\begin{align}\label{ineq:Lkp_abs}
\left| \int_{\Rd} \delta^{\mathfrak{K}_2} (t,x,y;z) \, \mathfrak{K}_1(z) J(z)dz
  \right| \leq c t^{-1} \rr_t(y-x)\,.
\end{align}
\end{lemma}
\pf
Let ${\rm I}$ be the left hand side of \eqref{ineq:Lkp_abs}.
We note that the integral defining ${\rm I}$ converges absolutely.
Using \eqref{eq:L_delta_r} with $r=h^{-1}(1/t)$, and \eqref{ineq:est_delta_1_crit},
\begin{align*}
{\rm I} &\leq  c\int_{|z|\geq  h^{-1}(1/t)} \aF_{2} (t,x,y;z) \,  \mathfrak{K}_1(z) J(z)dz
+c \int_{|z|< h^{-1}(1/t)} \aF_{1} (t,x,y;z) \,  \mathfrak{K}_1(z) J(z)dz\\
&\quad + \left| \int_{\Rd} z \left(\ind_{|z|<h^{-1}(1/t)} - \ind_{|z|<1}\right)   \mathfrak{K}_1(z) J(z)dz\right| |\nabla_x p^{\mathfrak{K}_2}(t,x,y)|\,.
\end{align*}
By \eqref{ineq:k-int_control} and Proposition~\ref{prop:gen_est_crit}, 
 the last term is bounded by $c\, t^{-1}\rr_t(y-x)$.
The same holds for the first two terms, because of \eqref{e:psi1}, \eqref{ineq:k-bounded}, 
\cite[Lemma~5.1 (8) and~5.9]{GS-2018}.
\qed

In what follows, 
we shall see a
difference in
the estimates
compared to \cite{GS-2018}.
The forthcoming
result is an analogue of
\cite[Theorem~2.9]{GS-2018}
suitable for the present development.

\begin{proposition}\label{thm:delta_crit}
Assume \eqref{e:psi1}, \eqref{eq:intro:wlsc}, \eqref{ineq:k-bounded}, \eqref{ineq:k-int_control}. For every $T>0$, the inequalities
\begin{align}
&\int_{\Rd} |\delta^{\mathfrak{K}} (t,x,y;z)|\, J(z)dz
\leq c\, \ct(t)\,  t^{-1} \rr_t(y-x)\,, \label{e:fract-der-est1-crit}\\
\int_{\Rd} |\delta^{\mathfrak{K}} (t,x',y;z)-&\delta^{\mathfrak{K}} (t,x,y;z)|\, J(z)dz \leq c \left(\frac{|x'-x|}{h^{-1}(1/t)} \land 1\right)  \ct(t)\, t^{-1} \big( \rr_t(y-x') + \rr_t(y-x)\big),\nonumber 
\end{align}
hold for all $t\in(0,T]$, $x,x',y\in\Rd$ with
\begin{enumerate}
\item[(a)] $\ct(t)=\TCh(t)$ and $c=c(d,T,\param)$ if $\lah=1$,
\item[(b)] $\ct(t)=t \,[h^{-1}(1/t)]^{-1}$ and $c=c(d,T,\param,\uah,c_h)$ if  \eqref{eq:intro:wusc}  holds for $0<\lah \leq \uah<1$.
\end{enumerate}
\end{proposition}
\pf
By \eqref{ineq:est_delta_1_crit} 
with $r=h^{-1}(1/t)$ we get
\begin{align*}
\int_{\Rd} &|\delta^{\mathfrak{K}} (t,x,y;z)|\, J(z)dz \\
&\leq \int_{\Rd} |\delta_{1.r}^{\mathfrak{K}} (t,x,y;z)|\, J(z)dz + \int_{\Rd} |z| \left| \ind_{|z|<r}-\ind_{|z|<1} \right| J(z)\, dz\, |\nabla_x p^{\mathfrak{K}}(t,x,y)| \\
&\leq c \int_{|z|\geq  h^{-1}(1/t)} \aF_{2} (t,x,y;z) \, J(z)dz
+c \int_{|z|< h^{-1}(1/t)} \aF_{1} (t,x,y;z) \, J(z)dz\\
&+\int_{\Rd} |z| \left| \ind_{|z|<h^{-1}(1/t)}-\ind_{|z|<1} \right| J(z) dz\, \left[h^{-1}(1/t)\right]^{-1} \rr_t(y-x)\,.
\end{align*}
The first inequality in the statement follows from \eqref{e:psi1},
\cite[Lemma~5.1 (8) and~5.9]{GS-2018} and Lemma~\ref{lem:int_J}.
Now we prove the second inequality.
If $|x'-x|\geq h^{-1}(1/t)$, then
\begin{align*}
\int_{\Rd} \left(|\delta^{\mathfrak{K}} (t,x',y;z)|+|\delta^{\mathfrak{K}} (t,x,y;z)|\right)J(z)dz
\leq c \, \ct(t) t^{-1} \left( \rr_t(y-x')+\rr_t(y-x) \right)\,.
\end{align*}
If $|x'-x|< h^{-1}(1/t)$, we use \eqref{ineq:diff_delta_1_crit}, \eqref{ineq:est_grad_1}
and, again, 
\cite[Lemma~5.1 and~5.9]{GS-2018} and Lemma~\ref{lem:int_J}.

\qed

The next result is a tool to relate the heat kernels corresponding to two coefficients $\mathfrak{K}_1$, $\mathfrak{K}_2$.
\begin{lemma}\label{lem:rozne_1}
Assume \eqref{e:psi1}, \eqref{eq:intro:wlsc} and let \eqref{ineq:k-bounded},
\eqref{ineq:k-int_control} hold for $\mathfrak{K}_1$, $\mathfrak{K}_2$, $\mathfrak{K}_3$.
For all $t>0$, $x,y\in\Rd$ and $s\in (0,t)$,
\begin{align*}
\frac{d}{d s} \int_{\Rd} &p^{\mathfrak{K}_1}(s,x,z) p^{\mathfrak{K}_2}(t-s,z,y)\,dz\\
&= \int_{\Rd} \LL_x^{\mathfrak{K}_1}p^{\mathfrak{K}_1}(s,x,z) \, p^{\mathfrak{K}_2}(t-s,z,y)\,dz
- \int_{\Rd} p^{\mathfrak{K}_1}(s,x,z)\, \LL_z^{\mathfrak{K}_2} p^{\mathfrak{K}_2}(t-s,z,y) \,dz\,,
\end{align*}
and
\begin{align*}
\int_{\Rd} \LL^{\mathfrak{K}_3}_x p^{\mathfrak{K}_1}(s,x,z)\,
p^{\mathfrak{K}_2}(t-s,z,y)\,dz
= &\int_{\Rd} p^{\mathfrak{K}_1}(s,x,z) \, \LL_z^{\mathfrak{K}_3} p^{\mathfrak{K}_2}(t-s,z,y)\,dz\,.
\end{align*}
\end{lemma}
\pf
The proof is the same as in \cite[Lemma~2.10]{GS-2018}. 
The first part follows by the dominated convergence theorem, which justifies the differentiation under the integral sign, and then by applying \eqref{eq:p_gen_klas}. The second identity is obtained after changing the order of integration and integrating by parts, see \eqref{e:delta-f-def} and \eqref{eq:L_delta}.
In both cases we use
the fact that for all $0<t_0<T<\infty$ there exists a constant $c=c(d,T,t_0,\param)$ such that for all $t\in[t_0,T]$, $x,y\in\Rd$,
\begin{align}\label{ineq:aux_Q0}
\int_{\Rd} |\delta^{\mathfrak{K}_1} (t,x,y;z)|\, J(z)dz
\leq c\, t^{-1}\rr_t(y-x)\leq c t_0^{-1}\rr_{t_0}(y-x)\,,
\end{align}
which is valid under the assumptions of the lemma, see \eqref{e:fract-der-est1-crit}.
\qed

\section{Analysis of the  heat kernel of $\LL^{\mathfrak{K}_w}$}
\label{sec:analysis_LL_2}

In this section we  work under $\Qzero$. In particular, we
consider $J(z)$  satisfying \eqref{e:psi1} and \eqref{eq:intro:wlsc}, and
$\kappa(x,z)$ such that \eqref{e:intro-kappa} and \eqref{e:intro-kappa-crit} hold.
For a fixed $w\in \Rd$ we let $\mathfrak{K}_w(z)=\kappa(w,z)$, as in \eqref{def:k-frozen}. Since \eqref{ineq:k-bounded} and \eqref{ineq:k-int_control} hold for $\mathfrak{K}_w$,  the results of
Section~\ref{sec:analysis_LL} remain in force.
Like in Section~\ref{sec:analysis_LL} we
let $p^{\mathfrak{K}_w}(t,x,y)$ be the heat kernel of  $\LL^{{\mathfrak K}_w}$.
This procedure is known as freezing coefficients of the operator $\LL^{\kappa}$ given in~\eqref{e:intro-operator-a1-crit1}.

First, we estimate $
\big( \LL_x^{\mathfrak{K}_{w'}}-\LL_x^{\mathfrak{K}_w} \big)p^{\mathfrak{K}}(t,x,y)
$. 
\begin{lemma}\label{lem:Lkp_abs-H}
Assume $\Qzero$ and let \eqref{ineq:k-bounded}, \eqref{ineq:k-int_control} hold for $\mathfrak{K}$. For every $T>0$ there exists a constant $c=c(d,T,\param,\kappa_2,\kb)$ such that
for all $t\in (0,T]$, $x,y, w, w'\in\Rd$ we have
\begin{align}\label{ineq:Lkp_abs-H}
\left| \int_{\Rd} \delta^{\mathfrak{K}} (t,x,y;z) \, \left( \kappa(w',z)- \kappa(w,z)\right) J(z)dz
  \right| \leq c \left( |w'-w|^{\beta}\land 1 \right) t^{-1} \rr_t(y-x)\,.
\end{align}
\end{lemma}
\pf
If $|w'-w|\geq 1$ we apply \eqref{ineq:Lkp_abs}.
Let ${\rm I}$ be the left hand side of \eqref{ineq:Lkp_abs-H}
and $|w'-w|< 1$.
We also note that the integral defining ${\rm I}$ converges absolutely.
Using 
\eqref{eq:L_delta_r} with $r=h^{-1}(1/t)$, 
and~\eqref{ineq:est_delta_1_crit},
\begin{align*}
{\rm I} &\leq  c\int_{|z|\geq  h^{-1}(1/t)} \aF_{2} (t,x,y;z) \, |\kappa(w',z)- \kappa(w,z)| J(z)dz\\
&\quad+c \int_{|z|< h^{-1}(1/t)} \aF_{1} (t,x,y;z) \, |\kappa(w',z)- \kappa(w,z)| J(z)dz\\
&\quad + \left| \int_{\Rd} z \left(\ind_{|z|<h^{-1}(1/t)} - \ind_{|z|<1}\right) \big( \kappa(w',z)- \kappa(w,z)\big) J(z)dz\right| |\nabla_x p^{\mathfrak{K}}(t,x,y)|\,.
\end{align*}
By
\eqref{e:intro-kappa-crit-H} and
 Proposition~\ref{prop:gen_est_crit} 
 the last term is bounded by $|w'-w|^{\beta}t^{-1}\rr_t(y-x)$.
The same is true for the first two terms by \eqref{e:psi1},
\eqref{e:intro-kappa-holder},  
\cite[Lemma~5.1 (8) and~5.9]{GS-2018}.
\qed

Now we estimate 
$\big(\LL_{x'}^{\mathfrak{K}_{w'}}-\LL_{x'}^{\mathfrak{K}_w} \big) p^{\mathfrak{K}}(t,x',y)-\big(\LL_{x}^{\mathfrak{K}_{w'}}-\LL_{x}^{\mathfrak{K}_w} \big) p^{\mathfrak{K}}(t,x,y)$.
\begin{lemma}\label{lem:Lkp_abs-H-H}
Assume $\Qzero$ and let \eqref{ineq:k-bounded}, \eqref{ineq:k-int_control} hold for $\mathfrak{K}$.
For every $T>0$ there exists a constant $c=c(d,T,\param,\kappa_2,\kb)$
such that for all $t\in(0,T]$, $x,x',y, w, w'\in\Rd$,
\begin{align}\label{ineq:Lkp_abs-H-H}
&\left| \int_{\Rd} \left( \delta^{\mathfrak{K}} (t,x',y;z)-\delta^{\mathfrak{K}} (t,x,y;z) \right) \left( \kappa(w',z)- \kappa(w,z)\right) J(z)dz  \right| \nonumber \\
 &\qquad \leq c \left(\frac{|x'-x|}{h^{-1}(1/t)} \land 1\right) \left( |w'-w|^{\beta}\land 1 \right) t^{-1} \big( \rr_t(y-x') + \rr_t(y-x)\big)\,.
\end{align}
\end{lemma}
\pf
If $|x'-x|\geq h^{-1}(1/t)$ we apply \eqref{ineq:Lkp_abs-H}. 
Let ${\rm I}$ be the left hand side of \eqref{ineq:Lkp_abs-H-H} and $|x'-x|< h^{-1}(1/t)$.
By \eqref{eq:L_delta_r}
with $r=h^{-1}(1/t)$, and \eqref{ineq:diff_delta_1_crit},
\begin{align*}
{\rm I} &\leq  c \left(\frac{|x'-x|}{h^{-1}(1/t)}\right) \int_{|z|\geq  h^{-1}(1/t)} \aF_{2} (t,x,y;z) \, |\kappa(w',z)- \kappa(w,z)| J(z)dz\\
&\quad+c  \left(\frac{|x'-x|}{h^{-1}(1/t)}\right) \int_{|z|< h^{-1}(1/t)} \aF_{1} (t,x,y;z) \, |\kappa(w',z)- \kappa(w,z)| J(z)dz\\
&\quad + \left| \int_{\Rd} z \left(\ind_{|z|<h^{-1}(1/t)} - \ind_{|z|<1}\right) \big( \kappa(w',z)- \kappa(w,z)\big) J(z)dz\right| |\nabla_{x'} p^{\mathfrak{K}}(t,x',y)-\nabla_{x}p^{\mathfrak{K}}(t,x,y)|\,.
\end{align*}
By \eqref{e:intro-kappa-crit-H}
and \eqref{ineq:est_grad_1},
 we bound the last expression by
$(|w'-w|^{\beta}\land 1 ) (|x'-x|/h^{-1}(1/t)) t^{-1}\rr_t(y-x)$.
For the first two terms we rely on \eqref{e:psi1},
\eqref{e:intro-kappa-holder},  
\cite[Lemma~5.1 (8) and~5.9]{GS-2018}.
\qed

In Lemma~\ref{lem:Lkp_abs-H}
and~\ref{lem:Lkp_abs-H-H}
our assumptions 
\eqref{ineq:k-bounded} and \eqref{ineq:k-int_control}
play an important role.
They also influence Proposition~\ref{prop:Hcont_kappa_crit1} and other results.
We note a difference in the estimates
\eqref{e:delta-difference-abs-crit1}
in comparison to the corresponding bound in
\cite[Theorem~2.11]{GS-2018}.

In what follows we provide several results on the regularity of the heat kernel $p^{\mathfrak{K}_{w}}(t,x,y)$ in~$w\in\Rd$.

\begin{proposition}\label{prop:Hcont_kappa_crit1}
Assume $\Qzero$. For every $T>0$ there exists a constant $c=c(d,T,\param,\kappa_2,\kb)$
such that
for all $t\in (0,T]$, $x,y,w,w'\in\Rd$,
\begin{align*}
|p^{\mathfrak{K}_{w'}}(t,x,y)-p^{\mathfrak{K}_w}(t,x,y)| &\leq c\, (|w'-w|^{\beta}\land 1)\,\rr_t(y-x)\,,\\ 
|\nabla_x p^{\mathfrak{K}_{w'}}(t,x,y)-\nabla_x p^{\mathfrak{K}_w}(t,x,y)| &\leq c  (|w'-w|^{\beta}\land 1) \left[h^{-1}(1/t)\right]^{-1} \rr_t(y-x) \,,\\ 
\left| \LL_x^{\mathfrak{K}_x} p^{\mathfrak{K}_{w'}}(t,x,y)- \LL_x^{\mathfrak{K}_x} p^{\mathfrak{K}_{w}}(t,x,y)\right|
&\leq c (|w'-w|^{\beta}\land 1)\,
 t^{-1}\rr_t(y-x) \,. 
\end{align*}
Moreover, for every $T>0$, the inequality
\begin{align}
\int_{\Rd} |\delta^{\mathfrak{K}_{w'}} (t,x,y;z)-\delta^{\mathfrak{K}_w} (t,x,y;z)| \,J(z)dz &\leq c (|w'-w|^{\beta}\land 1) \,\ct(t)\,
 t^{-1}\rr_t(y-x) \,, \label{e:delta-difference-abs-crit1}
\end{align}
holds for all $t\in (0,T]$, $x,y,w,w'\in\Rd$ with
\begin{enumerate}
\item[(a)] $\ct(t)=\TCh(t)$ and $c=c(d,T,\param,\kappa_2,\kb)$ if $\lah=1$,
\item[(b)] $\ct(t)=t \,[h^{-1}(1/t)]^{-1}$ and $c=c(d,T,\param,\kappa_2,\kb,\uah,c_h)$ if  \eqref{eq:intro:wusc}  holds for $0<\lah \leq \uah<1$.
\end{enumerate}
\end{proposition}
\pf
In what follows, we use 
\cite[Corollary~5.14, Lemma~5.6]{GS-2018}
and the monotonicity of $h^{-1}$ without further comment. The proof resembles that of \cite[Theorem~2.11]{GS-2018}, but in parts (ii),  (iii) and (iv) different adjustments are needed to use our assumptions. \\
(i) Using Lemma~\ref{lem:rozne_1}, we get
\begin{align*}
p^{\mathfrak{K}_{w'}}(t,x,y)-p^{\mathfrak{K}_{w}}(t,x,y)
&= 
\lim_{\varepsilon_1 \to 0^+} \int_{\varepsilon_1}^{t/2}
\int_{\Rd} p^{\mathfrak{K}_{w'}}(s,x,z) \left( \LL_z^{\mathfrak{K}_{w'}} 
-  \LL_z^{\mathfrak{K}_w}\right) p^{\mathfrak{K}_w}(t-s,z,y)\,dzds\\
&+
\lim_{\varepsilon_2\to 0^+ } \int_{t/2}^{t-\varepsilon_2}
\int_{\Rd} \left( \LL_x^{\mathfrak{K}_{w'}} 
-  \LL_x^{\mathfrak{K}_{w}}\right) p^{\mathfrak{K}_{w'}}(s,x,z)  p^{\mathfrak{K}_w}(t-s,z,y)\,dzds
\,.
\end{align*}
By Proposition~\ref{prop:gen_est_crit} and \eqref{ineq:Lkp_abs-H},
\begin{align*}
&\int_{\varepsilon}^{t/2}
\int_{\Rd}  p^{\mathfrak{K}_{w'}}(s,x,z)\,  | \!\left( \LL_z^{\mathfrak{K}_{w'}} 
-  \LL_z^{\mathfrak{K}_w}\right) p^{\mathfrak{K}_w}(t-s,z,y)|\,dzds\\
&\leq
c\, (|w'-w|^{\beta}\land 1) 
\int_{\varepsilon}^{t/2}
\int_{\Rd} \rr_s (z-x)\, (t-s)^{-1}\rr_{t-s}(y-z) \,dzds\\
&\leq c\, (|w'-w|^{\beta}\land 1)\,
\rr_t(y-x) \int_{\varepsilon}^{t/2} t^{-1}ds\,.
\end{align*}
Similarly,
\begin{align*}
&\int_{t/2}^{t-\varepsilon}
\int_{\Rd} |\!\left( \LL_x^{\mathfrak{K}_{w'}} 
-  \LL_x^{\mathfrak{K}_w}\right) p^{\mathfrak{K}_{w'}}(s,x,z) |\,  p^{\mathfrak{K}_w}(t-s,z,y)\,dzds
\leq c  \, (|w'-w|^{\beta}\land 1)\, \rr_t(y-x)\,.
\end{align*}

\noindent
(ii) 
Let $w_0\in\Rd$  be fixed. Define $\mathfrak{K}(z)=(\kappa_0/(2\kappa_1)) \kappa(w_0,z)$ and $\widehat{\mathfrak{K}}_w (z)=\mathfrak{K}_w(z)- \mathfrak{K}(z)$. By the construction of the L{\'e}vy process, we have
\begin{align}\label{eq:przez_k_0-impr}
p^{\mathfrak{K}_w}(t,x,y)=\int_{\Rd} p^{\mathfrak{K}}(t,x,\xi)
p^{\widehat{\mathfrak{K}}_w}(t,\xi,y)\,d\xi\,.
\end{align}
Then by  \eqref{ineq:est_diff_1} we can differentiate under the integral in \eqref{eq:przez_k_0-impr}. By Proposition~\ref{prop:gen_est_crit} we  get
\begin{align*}
|\nabla_x p^{\mathfrak{K}_{w'}}(t,x,y)-\nabla_x p^{\mathfrak{K}_w}(t,x,y)|
& \leq  \int_{\Rd} \left| \nabla_x p^{\mathfrak{K}}(t, x,\xi) \right|  \left| p^{\widehat{\mathfrak{K}}_{w'}}(t,\xi,y)-p^{\widehat{{\mathfrak{K}}}_w}(t,\xi,y)\right| d\xi\\
&\leq c (|w'-w|\land 1) \left[h^{-1}(1/t)\right]^{-1} \rr_t(y-x)\,.
\end{align*}

\noindent
(iii)
By \eqref{eq:przez_k_0-impr} we have
\begin{align*}
\delta^{\mathfrak{K}_{w'}} (t,x,y;z)-\delta^{\mathfrak{K}_w} (t,x,y;z)
=  \int_{\Rd} \delta^{\mathfrak{K}}(t,x,\xi;z) 
\left(p^{\widehat{\mathfrak{K}}_{w'}}(t,\xi,y)-p^{\widehat{{\mathfrak{K}}}_w}(t,\xi,y)\right)
 d\xi.
\end{align*}
Then by  \eqref{ineq:Lkp_abs},
\begin{align*}
\left| \LL_x^{\mathfrak{K}_x} p^{\mathfrak{K}_{w'}}(t,x,y)- \LL_x^{\mathfrak{K}_x} p^{\mathfrak{K}_{w}}(t,x,y)\right|
&\leq 
 \int_{\Rd} \left| \LL_x^{\mathfrak{K}_x}p^{\mathfrak{K}}(t,x,\xi ) \right|
\left| p^{\widehat{\mathfrak{K}}_{w'}}(t,\xi,y)-p^{\widehat{{\mathfrak{K}}}_w}(t,\xi,y)\right|
 d\xi\\
&\leq c (|w'-w|\land 1)\, t^{-1} \rr_{t}(y-x)\,.
\end{align*}

\noindent
(iv) By Proposition~\ref{thm:delta_crit},
\begin{align*}
\int_{\Rd} &|\delta^{\mathfrak{K}_{w'}} (t,x,y;z)-\delta^{\mathfrak{K}_w} (t,x,y;z)| \,J(z)dz \\
&\leq 
 \int_{\Rd} \left( \int_{\Rd}  |\delta^{\mathfrak{K}}(t,x,\xi;z)| \,J(z)dz\right)
\left| p^{\widehat{\mathfrak{K}}_{w'}}(t,\xi,y)-p^{\widehat{{\mathfrak{K}}}_w}(t,\xi,y)\right|
 d\xi\\
&\leq 
c (|w'-w|\land 1)\int_{\Rd} \ct(t) t^{-1}\rr_t(\xi-x)
\rr_t(y-\xi)\, d\xi \\
&\leq c(|w'-w|\land 1) \ct(t) t^{-1} \rr_t(y-x)\,.
\end{align*}
\qed

Our results mostly have the same form as those in \cite{GS-2018}, and similarly as in \cite{GS-2018} we are able
to deduce the joint continuity, the concentration of mass,
and cancellations.

\begin{lemma}\label{lem:cont_frcoef}
Assume $\Qzero$. The functions $p^{\mathfrak{K}_w}(t,x,y)$ and $\nabla_x p^{\mathfrak{K}_w}(t,x,y)$ 
are jointly continuous in $(t, x, y,w) \in (0,\infty)\times (\Rd)^3$.
The function $\LL_x^{\mathfrak{K}_{v}} p^{\mathfrak{K}_{w}}(t,x,y)$ is jointly continuous in $(t,x,y,w,v)\in (0,\infty)\times  (\Rd)^4$. Furthermore,
\begin{align}\label{e:some-estimates-2c-crit1}
\lim_{t \to 0^+ } \sup_{x\in\Rd} \left| \int_{\Rd} p^{\mathfrak{K}_y}(t,x,y)\, dy -1\right|=0
\end{align}
\end{lemma}
\pf
The result follows from Proposition~\ref{prop:Hcont_kappa_crit1},
\cite[Lemma~6.1]{GS-2018}, \eqref{ineq:est_delta_1_crit} and Lemma~\ref{l:convolution},
cf. \cite[Lemma~3.1, 3.2 and~3.4]{GS-2018}.
\qed

\begin{lemma}\label{e:some-estimates-2bb-crit1}
Assume $\Qzero$. Let $\beta_1\in [0,\beta]\cap [0,\lah)$.
For every $T>0$ there exists a constant $c=c(d,T,\param,\kappa_2,\kb,\beta_1)$
such that for all $t\in (0,T]$, $x\in\Rd$,
\begin{equation*}
 \left|\int_{\Rd} \nabla_x  p^{\mathfrak{K}_y} (t,x,y)\,dy \right|
\leq c\! \left[h^{-1}(1/t)\right]^{-1+\beta_1}\,.
\end{equation*}
\end{lemma}
\pf
The inequality stems from \eqref{ineq:est_diff_1},  Proposition~\ref{prop:Hcont_kappa_crit1} and Lemma~\ref{l:convolution},
cf. \cite[Lemma~3.4]{GS-2018}.

\qed

Compared to \cite{GS-2018},
the two expressions in
Lemma~\ref{l:some-estimates-3b-crit1} and Lemma~\ref{l:some-estimates-3b-crit1-impr}
need to be estimated separately.

\begin{lemma}\label{l:some-estimates-3b-crit1}
Assume $\Qzero$. Let $\beta_1\in [0,\beta]\cap [0,\lah)$.
For every $T>0$, the inequality
\begin{align*}
\int_{\Rd} \left|\int_{\Rd} \delta^{\mathfrak{K}_y} (t,x,y;z) \,dy \right| J(z)dz
&\leq c\, \ct(t)\, t^{-1}\left[h^{-1}(1/t)\right]^{\beta_1}, 
\end{align*}
holds for all $t\in (0,T]$, $x\in\Rd$ with
\begin{enumerate}
\item[(a)] $\ct(t)=\TCh(t)$ and $c=c(d,T,\param,\kappa_2,\kb,\beta_1)$ if $\lah=1$,
\item[(b)] $\ct(t)=t \,[h^{-1}(1/t)]^{-1}$ and $c=c(d,T,\param,\kappa_2,\kb,\beta_1,\uah,c_h)$ if  \eqref{eq:intro:wusc}  holds for $0<\lah \leq \uah<1$.
\end{enumerate}
\end{lemma}
\pf
We subtract zero and use
\eqref{e:delta-difference-abs-crit1}, to get
\begin{align*}
\int_{\Rd} \left|\int_{\Rd} \delta^{\mathfrak{K}_y} (t,x,y;z)
-\delta^{\mathfrak{K}_x} (t,x,y;z) \,dy \right| J(z)dz
\leq 
c \int_{\Rd} \ct(t) \err{0}{\beta_1}(t,x-y) \,dy\,.
\end{align*}
The result follows from Lemma~\ref{l:convolution}(a).
\qed

\begin{lemma}\label{l:some-estimates-3b-crit1-impr}
Assume $\Qzero$. Let $\beta_1\in [0,\beta]\cap [0,\lah)$.
For every $T>0$ there exists a constant $c=c(d,T,\param,\kappa_2,\kb,\beta_1)$ such that
for all $t\in (0,T]$, $x\in\Rd$,
\begin{align*}
\left| \int_{\Rd} \LL_x^{\mathfrak{K}_x} p^{\mathfrak{K}_y} (t,x,y) \,dy  \right|
&\leq c t^{-1}\left[h^{-1}(1/t)\right]^{\beta_1}\,.
\end{align*}
\end{lemma}
\pf
By Proposition~\ref{prop:Hcont_kappa_crit1}, we have
\begin{align*}
\left| \int_{\Rd} \LL_x^{\mathfrak{K}_x} p^{\mathfrak{K}_y} (t,x,y) \,dy  \right|=
\left| \int_{\Rd}\left( \LL_x^{\mathfrak{K}_x} p^{\mathfrak{K}_y} (t,x,y)- \LL_x^{\mathfrak{K}_x} p^{\mathfrak{K}_x} (t,x,y) \right)dy  \right|
\leq c \int_{\Rd} \err{0}{\beta_1}(t,x-y)\, dy\,.
\end{align*}
The result follows from Lemma~\ref{l:convolution}(a).
\qed

\section{Levi's construction of the heat kernel}\label{sec:constr}

In this section, we focus on the objects introduced in 
Section~\ref{sec:notation}
by the formulae
\eqref{e:q0-definition}--\eqref{e:def-phi-y-2}. The estimates are stated using  a short notation proposed in \eqref{def:err}.

\subsection{Construction of $q(t,x,y)$}\label{sec:q}

\begin{lemma}\label{l:estimates-q0-crit1} 
Assume $\Qzero$.
For every $T>0$ there exists a constant $c=c(d,T,\param,\kappa_2, \kb)\geq 1$ such that for all 
$\beta_1\in[0,\beta]$, $t\in (0,T]$ and $x,x',y,y'\in\Rd$
\begin{align}\label{e:q0-estimate-crit1}
|q_0(t,x,y)|\leq c \err{0}{\beta_1}(t,y-x)\,,
\end{align}
and for every $\gamma\in [0,\beta_1]$, 
\begin{align}
&|q_0(t,x,y)-q_0(t,x',y)|\nonumber\\
&\leq c \left(|x-x'|^{\beta_1-\gamma}\land 1\right)\left\{\left(\err{\gamma}{0}+\err{\gamma-\beta_1}{\beta_1}\right)(t,x-y)
+\left(\err{\gamma}{0}+\err{\gamma-\beta_1}{\beta_1}\right)(t,x'-y)\right\},\label{e:estimate-step3-crit1}
\end{align}
and
\begin{align}
&|q_0(t,x,y)-q_0(t,x,y')|\nonumber \\
&\leq c \left(|y-y'|^{\beta_1-\gamma}\land 1\right)\left\{\left(\err{\gamma}{0}+\err{\gamma-\beta_1}{\beta_1}\right)(t,x-y)
+\left(\err{\gamma}{0}+\err{\gamma-\beta_1}{\beta_1}\right)(t,x-y')\right\}.
\label{e:estimate-q0-2-crit1}
\end{align}
\end{lemma}
\pf
(i) \eqref{e:q0-estimate-crit1} follows from \eqref{ineq:Lkp_abs-H}.\\
(ii) For $|x-x'|\geq 1$ the inequality holds by \eqref{e:q0-estimate-crit1} and
\cite[(92)]{GS-2018}:
\begin{align*}
|q_0(t,x,y)|
\leq c  \err{0}{\beta_1}(t,y-x)
\leq c\left[ h^{-1}(1/T)\vee 1\right]^{\beta_1-\gamma} \err{\gamma-\beta_1}{\beta_1}(t,y-x)\,.
\end{align*}
For $1\geq |x-x'|\geq h^{-1}(1/t)$ the result follows from \eqref{e:q0-estimate-crit1} and
\begin{align*}
|q_0(t,x,y)|
\leq c \err{0}{\beta_1}(t,y-x) 
= c \left[ h^{-1}(1/t)\right]^{\beta_1-\gamma} \err{\gamma-\beta_1}{\beta_1}(t,y-x)
\leq c |x-x'|^{\beta_1-\gamma}  \err{\gamma-\beta_1}{\beta_1}(t,y-x)\,.
\end{align*}
Now,  \eqref{ineq:Lkp_abs-H} and \eqref{ineq:Lkp_abs-H-H}
yield
\begin{align*}
 & |q_0(t,x,y)-q_0(t,x',y)|=\left|\int_{\Rd} \delta^{\mathfrak{K}_y} (t,x,y;z)(\kappa(x,z)-\kappa(y,z))\,J(z)dz\right.\\
&  \hspace{0.1\linewidth} -\left. \int_{\Rd}\delta^{\mathfrak{K}_y}(t,x',y;z)(\kappa(x',z)-\kappa(y,z))\,J(z)dz\right|\\
& \hspace{0.05\linewidth} \leq  \left| \int_{\Rd}\left( \delta^{\mathfrak{K}_y}(t,x,y;z)-\delta^{\mathfrak{K}_y}(t,x',y;z)\right)  \left(\kappa(x,z)-\kappa(y,z)\right) J(z)dz\right| \\
& \hspace{0.1\linewidth} + \left| \int_{\Rd}\left( \delta^{\mathfrak{K}_y} (t,x',y;z)\right) \left(\kappa(x,z)-\kappa(x',z)\right) J(z)dz\right|\\
& \hspace{0.1\linewidth}  + c \left(|x-x'|^{\beta_1}\land 1\right)\int_{\Rd}|\delta^{\mathfrak{K}_y}(t,x',y;z)|\,J(z)dz\\
\leq c &\left(|x-y|^{\beta_1}\land 1\right) 
\left(\frac{|x-x'|}{h^{-1}(1/t)} \land 1\right) 
\big(\err{0}{0} (t,x-y)+\err{0}{0}(t,x'-y)\big)
+ c \left(|x-x'|^{\beta_1}\land 1\right) \err{0}{0}(t,x'-y).
\end{align*}
Applying 
$(|x-y|^{\beta_1}\land 1)\leq (|x-x'|^{\beta_1}\land 1) + (|x'-y|^{\beta_1}\land 1)$
we obtain
\begin{align*}
|q_0(t,x,y)-q_0(t,x',y)|\leq \ &c \left(\frac{|x-x'|}{h^{-1}(1/t)} \land 1\right)
\big(\err{0}{\beta_1} (t,x-y)+\err{0}{\beta_1}(t,x'-y)\big)\\
&+c \left(|x-x'|^{\beta_1}\land 1\right)\err{0}{0}(t,x'-y).
\end{align*}
Thus, in the last case $|x-x'|\leq  h^{-1}(1/t)\land 1$ we have
$|x-x'|/ h^{-1}(1/t)\leq |x-x'|^{\beta_1-\gamma} \left[h^{-1}(1/t)\right]^{\gamma-\beta_1}$
and $|x-x'|^{\beta_1}\leq |x-x'|^{\beta_1 -\gamma} \left[h^{-1}(1/t)\right]^{\gamma}$.\\
(iii) We treat the cases $|y-y'|\geq 1$ and $1\geq |y-y'|\geq h^{-1}(1/t)$ like in part (ii).
Note that by
$\delta^{\mathfrak{K}}(t,x,y;z)=\delta^{\mathfrak{K}}(t,-y,-x;z)$,
\eqref{ineq:Lkp_abs-H}, \eqref{ineq:Lkp_abs-H-H}
and Proposition~\ref{prop:Hcont_kappa_crit1},
\begin{align*}
&|q_0(t,x,y)-q_0(t,x,y')|\\
&\leq \left| \int_{\Rd} \delta^{\mathfrak{K}_y}(t,x,y;z)\left(\kappa(y',z)-\kappa(y,z)\right)J(z)dz\right| \\
& \ \ \ +\left| \int_{\Rd}\left(\delta^{\mathfrak{K}_y}(t,x,y;z)-\delta^{\mathfrak{K}_y}(t,x,y';z)\right)\left(\kappa(x,z)-\kappa(y',z)\right)J(z)dz \right|\\
&\ \ \ +\left|\int_{\Rd}\left(\delta^{\mathfrak{K}_y}(t,x,y';z)-\delta^{\mathfrak{K}_{y'}}(t,x,y';z)\right)\kappa(x,z) J(z)dz\right|\\
&\ \ \ +\left|-\int_{\Rd}\left(\delta^{\mathfrak{K}_y}(t,x,y';z)-\delta^{\mathfrak{K}_{y'}}(t,x,y';z)\right)\kappa(y',z) J(z)dz \right|\\
&\leq c \left( |y-y'|^{\beta_1}\land 1\right) \err{0}{0}(t,x-y)\\
&\quad +c  \left( |x-y'|^{\beta_1}\land 1\right)  \left(\frac{|y-y'|}{h^{-1}(1/t)} \land 1\right)   \left(\err{0}{0}(t,x-y)+\err{0}{0}(t,x-y')\right)\\
&\quad +  c \left( |y-y'|^{\beta_1}\land 1\right) \err{0}{0}(t,x-y') \,.
\end{align*}
Applying 
$(|x-y'|^{\beta_1}\land 1)\leq (|x-y|^{\beta_1}\land 1) + (|y-y'|^{\beta_1}\land 1)$
we obtain
\begin{align*}
|q_0(t,x,y)-q_0(t,x,y')|\leq \ & c
 \left(\frac{|y-y'|}{h^{-1}(1/t)} \land 1\right) 
\big(\err{0}{\beta_1} (t,x-y)+\err{0}{\beta_1}(t,x-y')\big)\\
&+c \left(|y-y'|^{\beta_1}\land 1\right) \big( \err{0}{0}(t,x-y)+\err{0}{0}(t,x-y')\big).
\end{align*}
This proves \eqref{e:estimate-q0-2-crit1} in the case $|y-y'|\leq  h^{-1}(1/t)\land 1$.
\qed

We thus estimated $q_0$. The estimates are of the same form as in \cite[Lemma~3.6]{GS-2018}.
Using Lemma~\ref{l:convolution} they propagate to
functions $q_n$ and $q$ defined in \eqref{e:qn-definition} and \eqref{def:q}, respectively, see the proof of Theorem~\ref{t:definition-of-q-crit1}.
We also note that $\beta_1$ is used in Lemma~\ref{l:estimates-q0-crit1} merely
for technical convenience, but becomes relevant when estimating $q_n$ and $q$.

We stress that
among others,
 the inequality \eqref{e:difference-q-estimate-crit1} plays a special role and is often used
to improve the integrability or bounds of singular functions, sometimes
 along with cancellations
like those proved in
Lemma~\ref{e:some-estimates-2bb-crit1}
or
Lemma~\ref{l:some-estimates-3b-crit1-impr},
see  comments ahead of
Lemma~\ref{lem:int_grad_phi} and Lemma~\ref{lem:some-est_gen_phi_xy-crit1}.

\begin{theorem}\label{t:definition-of-q-crit1}
Assume $\Qzero$.
The series in \eqref{def:q} 
is locally uniformly absolutely convergent on $(0, \infty)\times \R^d \times \R^d$ and solves the integral equation
\begin{align}\label{e:integral-equation-crit1}
q(t,x,y)=q_0(t,x,y)+\int_0^t \int_{\Rd}q_0(t-s,x,z)q(s,z,y)\, dzds\, .
\end{align}
Moreover, 
for every $T> 0$ and $\beta_1\in (0,\beta]\cap (0,\lah)$ there is a constant $c=c(d,T,\param,\kappa_2,\kb, \beta_1)$
 such that on $(0,T]\times\Rd\times\Rd$,
\begin{align}\label{e:q-estimate-crit1}
|q(t,x,y)|\leq c \big(\err{0}{\beta_1}+\err{\beta_1}{0}\big)(t,x-y)\,,
\end{align}
and for any $\gamma\in (0,\beta_1]$ there is $c=c(d,T,\param,\kappa_2,\kb, \beta_1,\gamma)$ such that
on $(0, T]\times \Rd \times \Rd$,
\begin{align}
&|q(t,x,y)-q(t,x',y)|\nonumber\\
&\leq c \left(|x-x'|^{\beta_1-\gamma}\land 1\right)
\left\{\big(\err{\gamma}{0}+\err{\gamma-\beta_1}{\beta_1}\big)(t,x-y)+\big(\err{\gamma}{0}+\err{\gamma-\beta_1}{\beta_1}\big)(t,x'-y)\right\}\,,
\label{e:difference-q-estimate-crit1}
\end{align}
and
\begin{align}
&|q(t,x,y)-q(t,x,y')|\nonumber\\
&\leq c \left(|y-y'|^{\beta_1-\gamma}\land 1\right)
\left\{\big(\err{\gamma}{0}+\err{\gamma-\beta_1}{\beta_1}\big)(t,x-y)+\big(\err{\gamma}{0}+\err{\gamma-\beta_1}{\beta_1}\big)(t,x-y')\right\}\,.
\label{e:difference-q-estimate_1-crit1}
\end{align}
\end{theorem}
\pf
The proof follows from Lemmas~\ref{l:estimates-q0-crit1}
and~\ref{l:convolution} -- it is the same as for \cite[Theorem~3.7]{GS-2018}.

\qed

\subsection{Properties of $\phi_y(t,x,s)$ and $\phi_y(t,x)$}\label{sec:phi}

We shall prove estimates for the integral part of~ \eqref{e:p-kappa}. We use the notation introduced in \eqref{e:phi-y-def} and \eqref{e:def-phi-y-2}.

\begin{lemma}\label{lem:phi_cont_xy-crit1}
Assume $\Qzero$. Let $\beta_1\in (0,\beta]\cap (0,\lah)$.
For every $T>0$ there exists a constant 
$c=c(d,T,\param,\kappa_2,\kb, \beta_1)$
such that for all $t\in (0,T]$, $x,y\in\Rd$,
\begin{align*}
|\phi_y(t,x)|\leq c t \big(\err{0}{\beta_1}+\err{\beta_1}{0}\big)(t,x-y)\,.
\end{align*}
For any $T>0$ and $\gamma \in [0,1]\cap [0,\lah)$ there exists a constant
$c=c(d,T,\param,\kappa_2,\kb, \beta_1,\gamma)$
 such that
for all $t\in (0,T]$, $x,x',y\in \Rd$,
\begin{align*}
|\phi_{y}(t,x)-\phi_{y}(t,x')|&\leq c (|x-x'|^{\gamma}\land 1) \, t \left\{ \big( \err{\beta_1-\gamma}{0}+\err{-\gamma}{\beta_1}\big)(t,x-y)+ \big( \err{\beta_1-\gamma}{0}+\err{-\gamma}{\beta_1}\big)(t,x'-y) \right\}.
\end{align*}
For any $T>0$ and $\gamma \in (0,\beta)$ there exists a constant
$c=c(d,T,\param,\kappa_2,\kb, \beta_1,\gamma)$
 such that
for all $t\in (0,T]$, $x,y,y'\in \Rd$,
\begin{align*}
|\phi_{y}(t,x)-\phi_{y'}(t,x)|&\leq c (|y-y'|^{\beta_1-\gamma}\land 1)\, t \left\{ \big( \err{\gamma}{0}+\err{\gamma-\beta_1}{\beta_1}\big)(t,x-y)+ \big( \err{\gamma}{0}+\err{\gamma-\beta_1}{\beta_1}\big)(t,x-y') \right\}.
\end{align*}
\end{lemma}
\pf
The proof follows from Lemma~\ref{lem:pkw_holder},
Proposition~\ref{prop:gen_est_crit},
Theorem~\ref{t:definition-of-q-crit1}
and
Lemma~\ref{l:convolution}, and is the same as in \cite[Lemma~3.8]{GS-2018}.
\qed

\begin{lemma}\label{lem:phi_cont_joint-crit1}
Assume $\Qzero$.
The function $\phi_y(t,x)$ is jointly continuous in $(t,x,y)\in (0,\infty)\times \Rd \times \Rd$.
\end{lemma}
\pf
The idea of the proof is the same as that of \cite[Lemma~3.9]{GS-2018} and
relies on
Proposition~\ref{prop:gen_est_crit},
\eqref{e:q-estimate-crit1},
\cite[(94)]{GS-2018},
Lemma~\ref{l:convolution} and~
\ref{lem:cont_frcoef},
\cite[Lemma~5.6 and~5.15]{GS-2018}.
\qed

From this moment on, the major effort is  to obtain sufficient regularity of the integral part of \eqref{e:p-kappa}.
Recall that
we need
$\beta_1<\lah$ in order to apply Lemma~\ref{l:convolution}. 
The additional condition $1<\beta_1+\lah$, known in certain contexts as the balance condition, which shall appear below in our assumptions,
makes it possible to differentiate \eqref{e:def-phi-y-2}, that is, to calculate and estimate its gradient, see Lemma~\ref{l:gradient-phi-y-crit1}. 
We need such a result if we want to apply either the strong or the weak operator \eqref{e:intro-operator-a1-crit1} to the candidate of a solution defined by \eqref{e:p-kappa}.

\begin{lemma}\label{lem:phi_pomoc-crit1}
Assume $\Qzero$. For all $0<s<t$, $x,y\in\Rd$,
\begin{align}\label{eq:grad_phi_pomoc-crit1}
\nabla_x \phi_y(t,x,s)=\int_{\Rd} \nabla_x p^{\mathfrak{K}_z}(t-s,x,z)q(s,z,y)\, dz\,,
\end{align}
\begin{align}\label{e:L-on-phi-y2-crit1}
\LL_x^{\mathfrak{K}_x}\phi_y(t,x,s)
=\int_{\Rd} \LL_x^{\mathfrak{K}_x}p^{\mathfrak{K}_z}(t-s,x,z)q(s,z,y)\, dz\,.
\end{align}
\end{lemma}
\pf
We get \eqref{eq:grad_phi_pomoc-crit1} by \eqref{ineq:est_diff_1}, \eqref{e:q-estimate-crit1},
Lemma~\ref{l:convolution}, and the dominated convergence theorem.
Now, by \eqref{e:phi-y-def} and
\eqref{eq:grad_phi_pomoc-crit1},
\begin{align}
\LL_x^{\mathfrak{K}_x}\phi_y(t,x,s)
=\int_{\Rd}  \left(\int_{\Rd} \delta^{\mathfrak{K}_z} (t-s, x,z;w) q(s,z, y)
\,dz\right) \kappa(x,w)J(w) dw\,. \label{e:L-on-phi-y2-first-crit1}
\end{align}
Finally, we use Fubini's theorem justified by \eqref{ineq:aux_Q0}, 
\eqref{e:q-estimate-crit1} and Lemma~\ref{l:convolution}(b).
\qed

In the proof of the next result we recognize a typical {\it modus operandi}
when dealing with integrals of functions that at first glance
 seem to be too singular:
we add and subtract $q(s,x,y)$,
use its regularity 
 \eqref{e:difference-q-estimate-crit1} 
(which reduces part of the singularity)
and profit from
cancellations, this time from Lemma~\ref{e:some-estimates-2bb-crit1}.

\begin{lemma}\label{lem:int_grad_phi}
Assume $\Qzero$ and $1-\lah <\beta\land \lah$.
Let $\beta_1\in (0,\beta]\cap (0,\lah)$.
For every $T>0$ 
  there exists a constant $c=c(d,T,\param,\kappa_2,\kb, \beta_1)$
such that for all $t\in (0,T]$, $x\in\Rd$,
\begin{align*}
\int_{\Rd} \int_0^t \left| \nabla_x \phi_y(t,x,s) \right| ds \, dy\leq c \left[h^{-1}(1/t)\right]^{-1+\beta_1}\,.
\end{align*}
\end{lemma}
\pf
By the monotonicity of $h^{-1}$, if we prove the statement for some value of $\beta_1$, then it also holds for smaller values.
We assume that $1-\lah<\beta_1$ and we let
 $\gamma \in (0,\beta_1)$ satisfying $1-\lah<\beta_1-\gamma$.
By 
\eqref{eq:grad_phi_pomoc-crit1},
Proposition~\ref{prop:gen_est_crit}, 
\eqref{e:difference-q-estimate-crit1}, Lemma~\ref{e:some-estimates-2bb-crit1}  and \eqref{e:q-estimate-crit1},
\begin{align*}
\left| \nabla_x \phi_y(t,x,s) \right|
 &\leq \int_{\Rd} \left| \nabla_x p^{\mathfrak{K}_z}(t-s,x,z) \right| \left| q(s,z,y) -q(s,x,y)\right| dz\\
&\quad + \left| \int_{\Rd} \nabla_x p^{\mathfrak{K}_z}(t-s,x,z)\, dz \right| \left| q(s,x,y)\right|\\
&\leq
\int_{\Rd}  (t-s)\err{-1}{\beta_1-\gamma}(t-s, x-z) \big(\err{\gamma}{0}+\err{\gamma-\beta_1}{\beta_1}\big)(s,z-y)\,dz \\
&\quad +\int_{\Rd}  (t-s)\err{-1}{\beta_1-\gamma}(t-s, x-z) \,dz\, \big(\err{\gamma}{0}+\err{\gamma-\beta_1}{\beta_1}\big)(s,x-y) \\
&\quad+  \left[ h^{-1}(1/(t-s))\right]^{-1+\beta_1} \big(\err{0}{\beta_1}+\err{\beta_1}{0}\big)(s,x-y)\,.
\end{align*}
Finally, we integrate 
in $y$ over $\Rd$
using Lemma~\ref{l:convolution}(a) and
then in $s$ over $(0,t)$ using \cite[Lemma~5.15]{GS-2018}.
Note that in the last step we integrate
$[h^{-1}(1/(t-s))]^{-1+\beta_1-\gamma}$, which requires a condition $(-1+\beta_1-\gamma)/\lah +1 >0$, equivalently $\lah+\beta_1>1+\gamma$, and is fulfilled thanks to our assumptions.
\qed

\begin{lemma}\label{l:gradient-phi-y-crit1}
Assume $\Qzero$ and $1-\lah <\beta\land \lah$.
For every $T>0$ 
there exists a constant 
$c=c(d,T,\param,\kappa_2,\kb,\beta)$
 such that
for all  $t \in(0,T]$, $x,y\in\Rd$,
\begin{align}\label{e:gradient-phi-y-crit1}
&\nabla_x\phi_y(t,x)=\int_0^t \int_{\Rd} \nabla_x p^{\mathfrak{K}_z}(t-s,x,z) q(s,z,y)\, dzds\,,\\
\nonumber\\
\label{e:gradient-phi-y-estimate-crit1}
&\left|\nabla_x\phi_y(t,x) \right|\leq c \!\left[ h^{-1}(1/t)\right]^{-1} t \,\err{0}{0}(t,x-y)\,.
\end{align}
\end{lemma}
\pf
The proof is like in \cite[Lemma~3.10]{GS-2018} and
rests on
\eqref{eq:grad_phi_pomoc-crit1},
Proposition~\ref{prop:gen_est_crit}, \eqref{e:q-estimate-crit1},
Lemma~\ref{l:convolution},
\cite[(93), (94), Lemma~5.3 and~5.15,
Proposition~5.8]{GS-2018}, \eqref{e:difference-q-estimate-crit1},
Lemma~\ref{e:some-estimates-2bb-crit1}, and the fact that $\lah>1/2$.
\qed

So far, 
in Lemma~\ref{l:gradient-phi-y-crit1}
we managed to calculate and estimate the gradient
of $\phi_y(t,x)$, which is
the integral part of \eqref{e:p-kappa}.
Now we shall treat in a similar fashion
the operator \eqref{e:intro-operator-a1-crit1} acting on
$\phi_y(t,x)$. 
The first step is to show that the operator can actually be applied (that the respective integrals converge) and to find a formula  --
Lemma~\ref{e:L-on-phi-y-crit1}. The second step is to prove the estimates in
Lemma~\ref{lem:I_0_oszagorne-crit1-impr}.
To achieve that, for the first step, in
Lemma~\ref{lem:some-est_gen_phi_xy-crit1}
and~\ref{ineq:I_0_oszagorne-crit1},
we prove auxiliary bounds 
justifying the use of Fubini's theorem,
however those technical results do not provide
the desired estimates for the second step.
The reason is that when using \eqref{e:difference-q-estimate-crit1},
due to the position of the absolute value, we cannot 
make use of additional cancellations
and we merely rely on
\eqref{e:fract-der-est1-crit} and
Lemma~\ref{l:some-estimates-3b-crit1},
which causes extra growth.
Therefore, contrary to \cite{GS-2018}, we are forced to distinguish between the two steps.
An improvement of the estimates, taking cancellations into account, is given in
Lemma~\ref{lem:some-est_gen_phi_xy-crit1-impr},
Corollary~\ref{cor:int_Lphi} and 
Lemma~\ref{lem:I_0_oszagorne-crit1-impr}.

In (a) below we address the critical case. In (b) we deal with the super-critical case.

\begin{lemma}\label{lem:some-est_gen_phi_xy-crit1}
Assume $\Qzero$. Let $\beta_1\in (0,\beta]\cap (0,\lah)$.
For all $T>0$, $\gamma \in(0,\beta_1]$ the inequalities
\begin{align}
\int_{\Rd}\left(\int_{\Rd} |\delta^{\mathfrak{K}_z} (t-s, x,z;w)||q(s,z, y)|
\,dz\right) \kappa(x,w)J(w) dw\nonumber \hspace{0.15\linewidth}\\  
 \leq c_1\int_{\Rd} \ct(t-s)\err{0}{0}(t-s, x-z)
\big(\err{0}{\beta_1}+\err{\beta_1}{0}\big)(s,z-y)\,dz\,,   \label{e:Fubini1-crit1} \\
\nonumber \\
\int_{\Rd}  \left| \int_{\Rd} \delta^{\mathfrak{K}_z}(t-s,x,z;w)q(s,z,y)\,dz \right| \kappa(x,w)J(w)dw
\leq c_2 \big( {\rm I}_1+{\rm I}_2+{\rm I}_3 \big), \label{ineq:some-est_gen_phi_xy-crit1}
\end{align}
where 
\begin{align*}
{\rm I}_1+{\rm I}_2+{\rm I}_3:=
& \int_{\Rd} \ct(t-s) \err{0}{\beta_1-\gamma}(t-s,x-z)  \big(\err{\gamma}{0}+\err{\gamma-\beta_1}{\beta_1}\big)(s,z-y) \,dz \\
& + \ct(t-s)
(t-s)^{-1} \left[h^{-1}(1/(t-s))\right]^{\beta_1-\gamma} 
\big(\err{\gamma}{0}+\err{\gamma-\beta_1}{\beta_1}\big)(s,x-y) \\
& +  \,\ct(t-s) (t-s)^{-1}\left[h^{-1}(1/(t-s))\right]^{\beta_1}  \big(\err{0}{\beta_1}+\err{\beta_1}{0}\big)(s,x-y)\,,
\end{align*}
hold for all  $0<s<t\leq T$, $x,y\in\Rd$ with
\begin{enumerate}
\item[(a)] $\ct(t)=\TCh(t)$ and $c_1=c_1(d,T,\param,\kappa_2,\kb,\beta_1)$,  $c_2=c_2(d,T,\param,\kappa_2,\kb,\beta_1,\gamma)$ if $\lah=1$,
\item[(b)] $\ct(t)=t \,[h^{-1}(1/t)]^{-1}$ and $c_1=c_1(d,T,\param,\kappa_2,\kb,\beta_1,\uah,c_h)$, $c_2=c_2(d,T,\param,\kappa_2,\kb,\beta_1,\gamma,\uah,c_h)$ if  \eqref{eq:intro:wusc}  holds for $0<\lah \leq \uah<1$.
\end{enumerate}
\end{lemma}
\pf
The inequality \eqref{e:Fubini1-crit1} follows from \eqref{e:intro-kappa}, \eqref{e:fract-der-est1-crit} and \eqref{e:q-estimate-crit1}.
Next, 
let ${\rm I}_0$ be the left hand side of \eqref{ineq:some-est_gen_phi_xy-crit1}.
By \eqref{e:difference-q-estimate-crit1}, \eqref{e:q-estimate-crit1}, \eqref{e:intro-kappa},  \eqref{e:fract-der-est1-crit},
Lemma~\ref{l:some-estimates-3b-crit1} and~\ref{l:convolution}(a),
\begin{align*}
{\rm I}_0&\leq \int_{\Rd}  \int_{\Rd} |\delta^{\mathfrak{K}_z}(t-s,x,z;w)| |q(s,z,y)-q(s,x,y)|\,dz\, \kappa(x,w)J(w)dw\\
&\quad + \int_{\Rd}  \left| \int_{\Rd} \delta^{\mathfrak{K}_z}(t-s,x,z;w) \,dz\right| \kappa(x,w)J(w)dw \, |q(s,x,y)|\\
&\leq c \int_{\Rd} \left(  \int_{\Rd} |\delta^{\mathfrak{K}_z}(t-s,x,z;w)|\,J(w)dw \right) \left(|x-z|^{\beta_1-\gamma}\land 1\right) \big(\err{\gamma}{0}+\err{\gamma-\beta_1}{\beta_1}\big)(s,z-y) \,dz \\
&\quad + c \int_{\Rd} \left(  \int_{\Rd} |\delta^{\mathfrak{K}_z}(t-s,x,z;w)|\,J(w)dw \right) \left(|x-z|^{\beta_1-\gamma}\land 1\right) dz\, \big(\err{\gamma}{0}+\err{\gamma-\beta_1}{\beta_1}\big)(s,x-y) \\
&\quad + c  \, \ct(t-s)\, (t-s)^{-1}\left[h^{-1}(1/(t-s))\right]^{\beta_1}  \big(\err{0}{\beta_1}+\err{\beta_1}{0}\big)(s,x-y)
\leq c ({\rm I}_1+{\rm I}_2+{\rm I}_3)\,.
\end{align*}
\qed

The inequality 
\eqref{ineq:some-est_gen_phi_xy-crit1}
looks a bit rough, but it is left in such form on purpose: it is used not only to prove Lemma~\ref{ineq:I_0_oszagorne-crit1} and Lemma~\ref{lem:some-est_p_kappa-crit1}, but also to
 shorten the argument in the proof of 
Lemma~\ref{lem:I_0_oszagorne-crit1-impr}.

\begin{lemma} \label{ineq:I_0_oszagorne-crit1}
Assume $\Qzero$ and $1-\lah <\beta\land \lah$.
For any $\beta_1\in (0,\beta]$ such that $1-\lah<\beta_1<\lah$ and $0<\gamma_1\leq \gamma_2\leq \beta_1$ satisfying
$$
1-\lah<\beta_1-\gamma_1\,,\qquad\qquad 2\beta_1-\gamma_2<\lah\,,
$$
the inequality
\begin{align*}
\int_{\Rd}  \int_0^t &\left| \int_{\Rd} \delta^{\mathfrak{K}_z}(t-s,x,z;w)q(s,z,y)\,dz \right| ds\,\kappa(x,w)J(w)dw\nonumber \\
&\hspace{0.3\linewidth}\leq c  \,\ct(t) \big(\err{0}{\beta_1}+\err{\gamma_1}{\beta_1-\gamma_1}+\err{\beta_1+\gamma_1-\gamma_2}{0}\big)(t,x-y)\,,
\end{align*}
holds for all  $t\in(0,T]$, $x,y\in\Rd$ with
\begin{enumerate}
\item[(a)] $\ct(t)=\TCh(t)$ and $c=c(d,T,\param,\kappa_2,\kb,\beta_1,\gamma_1,\gamma_2)$ if $\lah=1$,
\item[(b)] $\ct(t)=t \,[h^{-1}(1/t)]^{-1}$ and $c=c(d,T,\param,\kappa_2,\kb,\beta_1,\gamma_1,\gamma_2,\uah,c_h)$ if  \eqref{eq:intro:wusc}  holds for $0<\lah \leq \uah<1$.
\end{enumerate}
\end{lemma}
\pf
Let ${\rm I}_0$ be the left hand side of \eqref{ineq:some-est_gen_phi_xy-crit1}.
In the  two cases discussed below, we apply
Lemma~\ref{l:convolution}(b),  the monotonicity of $h^{-1}$ 
and $\TCh$ (see also Lemma~\ref{lem:cal_TCh}),
and $\Ab$ of 
\cite[Lemma~5.3]{GS-2018}.
For $s\in (0,t/2]$ we use \eqref{e:Fubini1-crit1}
to get
\begin{align*}
{\rm I}_0& 
\leq c\, \ct(t-s)\bigg\{\left( (t-s)^{-1}\left[h^{-1}(1/(t-s))\right]^{\beta_1}+(t-s)^{-1}\left[h^{-1}(1/s)\right]^{\beta_1}
+s^{-1}\left[h^{-1}(1/s)\right]^{\beta_1} \right) \\
&\hspace{0.52\linewidth} \times \,\err{0}{0}(t,x-y) + (t-s)^{-1}\err{0}{\beta_1}(t,x-y) \bigg\}\\
&\leq c\,\ct(t) \bigg\{\left( t^{-1}\left[h^{-1}(1/t)\right]^{\beta_1} +s^{-1} \left[h^{-1}(1/s)\right]^{\beta_1} \right) \err{0}{0}(t,x-y) + t^{-1}\err{0}{\beta_1}(t,x-y) \bigg\}.
\end{align*}
For $s\in (t/2,t)$ we use  \eqref{ineq:some-est_gen_phi_xy-crit1} with $\gamma=\gamma_1$.
While estimating the expression
\begin{align*}
\ct(t-s) \int_{\Rd} \err{0}{\beta_1-\gamma}(t-s,x-z)  \err{\gamma-\beta_1}{\beta_1}(s,z-y) \,dz\,,
\end{align*}
we use Lemma~\ref{l:convolution}(b) with $n_1=n_2=2\beta_1-\gamma_2$, $m_1=\beta_1-\gamma_1$, $m_2=\beta_1$ and 
later Lemma~\ref{lem:cal_TCh},
so our assumptions concerning the choice of $\gamma_1$, $\gamma_2$ are used.
More precisely, we have
\begin{align*}
{\rm I}_1 &
\leq  c \ct(t-s)\bigg\{ (t-s)^{-1} \left[h^{-1}(1/(t-s))\right]^{\beta_1-\gamma_1}  \left[h^{-1}(1/s)\right]^{\gamma_1} + s^{-1} \left[h^{-1}(1/s)\right]^{\beta_1+\gamma_1-\gamma_2}  \\
&\hspace{0.32\linewidth}+  (t-s)^{-1} \left[h^{-1}(1/(t-s))\right]^{2\beta_1-\gamma_2}  \left[h^{-1}(1/s)\right]^{\gamma_1-\beta_1}  \bigg\}  \err{0}{0}(t,x-y)\\
&\quad +  c \ct(t-s)   (t-s)^{-1} \left[h^{-1}(1/(t-s))\right]^{\beta_1-\gamma_1}  \left[h^{-1}(1/s)\right]^{\gamma_1-\beta_1} \err{0}{\beta_1}(t,x-y)\\
&\quad +c \ct(t-s) s^{-1} \left[h^{-1}(1/s)\right]^{\gamma_1} \err{0}{\beta_1-\gamma_1}(t,x-y) \\
& \leq  c \ct(t-s) \bigg\{ (t-s)^{-1} \left[h^{-1}(1/(t-s))\right]^{\beta_1-\gamma_1}  \left[h^{-1}(1/t)\right]^{\gamma_1} + t^{-1} \left[h^{-1}(1/t)\right]^{\beta_1+\gamma_1-\gamma_2}\\
&\hspace{0.32\linewidth}+(t-s)^{-1} \left[h^{-1}(1/(t-s))\right]^{2\beta_1-\gamma_2}  \left[h^{-1}(1/t)\right]^{\gamma_1-\beta_1}  \bigg\}  \err{0}{0}(t,x-y)\\
&\quad + c \ct(t-s)   (t-s)^{-1} \left[h^{-1}(1/(t-s))\right]^{\beta_1-\gamma_1}  \left[h^{-1}(1/t)\right]^{\gamma_1-\beta_1} \err{0}{\beta_1}(t,x-y)\\
&\quad +c \ct(t-s) t^{-1} \left[h^{-1}(1/t)\right]^{\gamma_1} \err{0}{\beta_1-\gamma_1}(t,x-y)\,.
\end{align*}
Next, 
like above with \cite[(94)]{GS-2018},
\begin{align*}
{\rm I}_2 &
\leq c \,   \ct(t-s) (t-s)^{-1} \left[h^{-1}(1/(t-s))\right]^{\beta_1-\gamma_1} \big(\err{\gamma_1}{0}+\err{\gamma_1-\beta_1}{\beta_1}\big)(t,x-y)\,.
\end{align*}
Similarly,
${\rm I}_3\leq c \,\ct(t-s) (t-s)^{-1} \left[h^{-1}(1/(t-s))\right]^{\beta_1}  \big(\err{0}{\beta_1}+\err{\beta_1}{0}\big)(t,x-y)$.
Finally, by
\cite[Lemma~5.15]{GS-2018}
and Lemma~\ref{lem:cal_TCh},
and a fact that $\lah>1/2$,
\begin{align*}
\int_0^t {\rm I}_0\,ds 
\leq c \,\ct(t) \big(\err{0}{\beta_1}+\err{\gamma_1}{\beta_1-\gamma_1}+\err{\beta_1+\gamma_1-\gamma_2}{0}\big)(t,x-y)\,.
\end{align*}
\qed

We can now successfully apply
\eqref{e:intro-operator-a1-crit1} to
\eqref{e:def-phi-y-2}.

\begin{lemma}\label{e:L-on-phi-y-crit1}
Assume $\Qa$ or $\Qb$.
We have for all  $t >0$, $x,y\in\Rd$,
\begin{equation*}
\LL_x^{\mathfrak{K}_x} \phi_y(t,x)= \int_0^t \int_{\Rd} \LL_x^{\mathfrak{K}_x} p^{\mathfrak{K}_z}(t-s,x,z) q(s,z,y)\, dzds\,.
\end{equation*}
\end{lemma}
\pf
By \eqref{e:def-phi-y-2} and \eqref{e:gradient-phi-y-crit1}
in the first equality, and  Lemma~\ref{ineq:I_0_oszagorne-crit1} and \eqref{e:Fubini1-crit1} in the second
(allowing us to change the order of integration twice) the proof is as follows
\begin{align*}
\LL_x^{\mathfrak{K}_x} \phi_y(t,x)
&=\int_{\Rd} \left( \int_0^t \int_{\Rd} \delta^{\mathfrak{K}_z}(t-s,x,z;w)q(s,z,y)\,dzds\right) \kappa(x,w)J(w)dw\\
&=  \int_0^t \int_{\Rd}\left( \int_{\Rd} \delta^{\mathfrak{K}_z}(t-s,x,z;w)\, \kappa(x,w)J(w)dw\right) q(s,z,y)\,dzds\,.
\end{align*}
\qed

We improve the estimates.

\begin{lemma}\label{lem:some-est_gen_phi_xy-crit1-impr}
Assume $\Qzero$.
Let $\beta_1\in (0,\beta]\cap (0,\lah)$.
For all $T>0$, $\gamma \in(0,\beta_1]$ 
there exist  constants $c_1=c_1(d,T,\param,\kappa_2,\kb,\beta_1)$ and $c_2=c_2(d,T,\param,\kappa_2,\kb,\beta_1,\gamma)$ such that
for all  $0<s<t\leq T$, $x,y\in\Rd$,
\begin{align}
\left| \LL_x^{\mathfrak{K}_x}\phi_y(t,x,s)\right|
& \leq c_1\int_{\Rd} \err{0}{0}(t-s, x-z) \big(\err{0}{\beta_1}+\err{\beta_1}{0}\big)(s,z-y)\,dz\,,   \label{e:Fubini1-crit1-impr} \\
\nonumber \\
\left| \LL_x^{\mathfrak{K}_x}\phi_y(t,x,s)\right|
&\leq c_2 \big( {\rm I}_1+{\rm I}_2+{\rm I}_3 \big), \label{ineq:some-est_gen_phi_xy-crit1-impr}
\end{align}
where 
\begin{align*}
{\rm I}_1+{\rm I}_2+{\rm I}_3:=
& \int_{\Rd}  \err{0}{\beta_1-\gamma}(t-s,x-z) \big(\err{\gamma}{0}+\err{\gamma-\beta_1}{\beta_1}\big)(s,z-y) \,dz \\
& +
(t-s)^{-1} \left[h^{-1}(1/(t-s))\right]^{\beta_1-\gamma}
\big(\err{\gamma}{0}+\err{\gamma-\beta_1}{\beta_1}\big)(s,x-y) \\
& +  \, (t-s)^{-1}\left[h^{-1}(1/(t-s))\right]^{\beta_1}  \big(\err{0}{\beta_1}+\err{\beta_1}{0}\big)(s,x-y)\,.
\end{align*}
\end{lemma}
\pf
The first inequality follows from \eqref{ineq:Lkp_abs} and \eqref{e:q-estimate-crit1}.
By \eqref{e:L-on-phi-y2-crit1},
\eqref{e:difference-q-estimate-crit1}, \eqref{e:q-estimate-crit1},
\eqref{ineq:Lkp_abs},
Lemma~\ref{l:some-estimates-3b-crit1-impr} and~\ref{l:convolution}(a),
\begin{align*}
\left| \LL_x^{\mathfrak{K}_x}\phi_y(t,x,s)\right|&\leq \int_{\Rd}  
\left| \LL_x^{\mathfrak{K}_x}p^{\mathfrak{K}_z}(t-s,x,z)\right|
 |q(s,z,y)-q(s,x,y)|\,dz \\
&\quad +   \left| \int_{\Rd}\LL_x^{\mathfrak{K}_x}p^{\mathfrak{K}_z}(t-s,x,z) \,dz\right| \, |q(s,x,y)|\\
&\leq c \int_{\Rd} \err{0}{0}(t-s,x-z) \left(|x-z|^{\beta_1-\gamma}\land 1\right) \big(\err{\gamma}{0}+\err{\gamma-\beta_1}{\beta_1}\big)(s,z-y) \,dz \\
&\quad + c \int_{\Rd}  \err{0}{0}(t-s,x-z)  \left(|x-z|^{\beta_1-\gamma}\land 1\right) dz\, \big(\err{\gamma}{0}+\err{\gamma-\beta_1}{\beta_1}\big)(s,x-y) \\
&\quad + c  \, (t-s)^{-1}\left[h^{-1}(1/(t-s))\right]^{\beta_1}  \big(\err{0}{\beta_1}+\err{\beta_1}{0}\big)(s,x-y)
\leq c ({\rm I}_1+{\rm I}_2+{\rm I}_3)\,.
\end{align*}
\qed

Here is a consequence of \eqref{e:Fubini1-crit1-impr},  \eqref{ineq:some-est_gen_phi_xy-crit1-impr}, Lemma~\ref{l:convolution}(a)
and \cite[Lemma~5.15]{GS-2018}.
\begin{corollary}\label{cor:int_Lphi}
Assume $\Qzero$ and $1-\lah <\beta\land \lah$.
Let $\beta_1 \in (0,\beta]\cap (0,\lah)$.
 For every $T>0$ there exists a constant 
$c=c(d,T,\param,\kappa_2,\kb, \beta_1)$
 such that for all $t\in(0,T]$, $x\in\Rd$,
\begin{align*}
\int_{\Rd} \int_0^t    \left|  \LL_x^{\mathfrak{K}_x}\phi_y(t,x,s)  \right|   ds \,dy \leq c  t^{-1} \left[h^{-1}(1/t)\right]^{\beta_1} \,.
\end{align*}
\end{corollary}

\begin{lemma}\label{lem:I_0_oszagorne-crit1-impr}
Assume $\Qzero$.
Let $\beta_1\in (0,\beta]\cap (0,\lah)$.
For all $T>0$, $0<\gamma_1\leq \gamma_2\leq \beta_1$ satisfying
$$
0<\beta_1-\gamma_1\,,\quad \qquad 2\beta_1-\gamma_2<\lah\,,
$$
there exists a constant $c=c(d,T,\param,\kappa_2,\kb,\beta_1,\gamma_1,\gamma_2)$ such that for all $t\in(0,T]$, $x,y\in\Rd$,
\begin{align}
\int_0^t \left| \LL_x^{\mathfrak{K}_x}\phi_y(t,x,s)\right|ds
\leq c \big(\err{0}{\beta_1}+\err{\gamma_1}{\beta_1-\gamma_1}+\err{\beta_1+\gamma_1-\gamma_2}{0}\big)(t,x-y)\,. \label{ineq:I_0_oszagorne-crit1-impr}
\end{align}
\end{lemma}
\pf
The proof goes by the same lines as the proof of Lemma~\ref{ineq:I_0_oszagorne-crit1} but with $\ct$ replaced by $1$,
and Lemma~\ref{lem:some-est_gen_phi_xy-crit1-impr} in place of 
Lemma~\ref{lem:some-est_gen_phi_xy-crit1}.
\qed
\begin{lemma}\label{lem:Lphi_cont-crit}
Assume $\Qa$ or $\Qb$.
The function 
$\LL_x^{\mathfrak{K}_x} \phi_y(t,x)$  is jointly continuous in   $(t,x,y)\in (0,\infty)\times \Rd\times \Rd$.
\end{lemma}
\pf
The proof is the same as in \cite[Lemma~3.13]{GS-2018} and requires
Lemma~\ref{e:L-on-phi-y-crit1}, \eqref{e:L-on-phi-y2-crit1}, \eqref{ineq:some-est_gen_phi_xy-crit1-impr}, \eqref{ineq:Lkp_abs},
\eqref{e:difference-q-estimate_1-crit1}, \eqref{e:q-estimate-crit1}, 
\cite[(94), Lemma~5.15]{GS-2018}, Lemmas~\ref{l:convolution} and~\ref{lem:cont_frcoef}.
\qed

In the final results of that section, we prepare
to calculate the time derivative of~\eqref{e:p-kappa}.

\begin{proposition}\label{l:phi-y-abs-cont-crit1}
Assume $\Qzero$.
For all $t>0$, $x,y\in \Rd$, $x\neq y$, we have 
\begin{align*}
\phi_y(t,x) 
=\int_0^t \left(q(r,x,y)+  \int_0^r \int_{\Rd} \LL_x^{\mathfrak{K}_z} p^{\mathfrak{K}_z}(r-s,x,z) q(s,z,y)\, dzds\right) dr\,.
\end{align*}
\end{proposition}
\pf
The idea of the proof is that
 differentiating
\eqref{e:def-phi-y-2} in $t>0$, we expect 
to get
\begin{align*}
\partial_t \phi_y(t,x) = q(t,x,y)+
 \int_0^t \int_{\Rd}\partial_t \, p^{\mathfrak{K}_z}(t-s,x,z) q(s,z,y)\, dzds\,.
\end{align*}
Actually, we intend to prove the following integral counterpart,
\begin{align*}
\phi_y(t,x) 
=\int_0^t \left(q(r,x,y)+  \int_0^r \int_{\Rd} \partial_r\, p^{\mathfrak{K}_z}(r-s,x,z) q(s,z,y)\, dzds\right) dr\,,
\end{align*}
see \eqref{eq:p_gen_klas}.
Therefore, the aim is to justify
\begin{align*}
\int_0^t & \int_0^r \int_{\Rd} \partial_r\, p^{\mathfrak{K}_z}(r-s,x,z) q(s,z,y)\, dzds\,dr=
\int_0^t \int_0^r \partial_r
\phi_y(r ,x,s)\, dsdr\\
&= 
\int_0^t \int_s^t \partial_r
\phi_y(r ,x,s)\, drds
= \int_0^t \left( \phi_y(t ,x,s)-
\lim_{\varepsilon \to 0^+} \phi_y(s+\varepsilon ,x,s) \right)ds\\
&= \phi_y(t,x)  - \int_0^t \lim_{\varepsilon \to 0^+} \phi_y(s+\varepsilon ,x,s)ds\,,
\end{align*}
and prove that 
$\lim_{\varepsilon \to 0^+} \phi_y(s+\varepsilon ,x,s) = q(s,x,y)$.
Details are
like in the proof of \cite[Lemma~3.14]{GS-2018}: we use
\eqref{e:phi-y-def}, \eqref{eq:p_gen_klas}, \eqref{ineq:Lkp_abs}, \eqref{e:q-estimate-crit1}, \eqref{e:L-on-phi-y2-crit1}, \eqref{ineq:I_0_oszagorne-crit1-impr}, \eqref{e:q0-estimate-crit1},
Lemma~\ref{l:convolution},
\cite[(92), (93), (94)]{GS-2018}, \eqref{e:some-estimates-2c-crit1},
\eqref{e:estimate-step3-crit1}, Proposition~\ref{prop:gen_est_crit},
\cite[Lemma~5.6]{GS-2018},
\cite[Theorem~7.21]{MR924157}.
\qed

If we combine
\eqref{e:integral-equation-crit1}
and
Lemma~\ref{e:L-on-phi-y-crit1}, then we can  represent the integrand in the statement of
Proposition~\ref{l:phi-y-abs-cont-crit1}
as $q_0(s,x,y)+ \LL_x^{\mathfrak{K}_x} \phi_y (s,x)$. Hence we conclude what follows.

\begin{corollary}\label{e:phi-y-partial_1-crit}
Assume $\Qa$ or $\Qb$.
For all $x,y\in\Rd$, $x\neq y$, the function $\phi_y(t,x)$ is differentiable in $t>0$  and
\begin{align*}
 \partial_t \phi_y(t,x)  =
q_0(t,x,y)+  \LL_x^{\mathfrak{K}_x} \phi_y (t,x)\,.
\end{align*}
\end{corollary}

\subsection{Properties of $p^\kappa(t, x, y)$}\label{sec:p_kappa}

We collect what can already be said about $p^\kappa(t, x, y)$.

\begin{lemma}\label{lem:some-est_p_kappa-crit1}
Assume $\Qzero$ and $1-\lah <\beta\land \lah$. Let $\beta_1\in (0,\beta]\cap (0,\lah)$.
For every $T>0$ 
the inequalities
\begin{align}
\int_{\Rd} | \delta^{\kappa}(t,x,y;z) | \,\kappa(x,z)J(z)dz &\leq c_1\, \ct(t) \err{0}{0}(t,x-y)\,,
\label{ineq:some-est_p_kappa-crit1} \\
\int_{\Rd} \left|\int_{\Rd} \delta^{\kappa}(t,x,y;z)\, dy \right|\kappa(x,z)&J(z)dz \leq c_2\,  \ct(t) t^{-1} \left[h^{-1}(1/t)\right]^{\beta_1}\,, \label{ineq:some-est_p_kappa_1-crit1}
\end{align}
hold for all  $t\in(0,T]$, $x,y\in\Rd$ with
\begin{enumerate}
\item[(a)] $\ct(t)=\TCh(t)$ and $c_1=c_1(d,T,\param,\kappa_2,\kb,\beta)$, $c_2=c_2(d,T,\param,\kappa_2,\kb,\beta_1)$ if $\lah=1$,
\item[(b)] $\ct(t)=t \,[h^{-1}(1/t)]^{-1}$ and $c_1=c_1(d,T,\param,\kappa_2,\kb,\beta,\uah,c_h)$, $c_2=c_2(d,T,\param,\kappa_2,\kb,\beta_1,\uah,c_h)$ if  \eqref{eq:intro:wusc}  holds for $0<\lah \leq \uah<1$.
\end{enumerate}
\end{lemma}
\pf
By \eqref{e:p-kappa} and \eqref{e:gradient-phi-y-crit1},
\begin{align*}
\delta^{\kappa}(t,x,y;w)=\delta^{\mathfrak{K}_y}(t,x,y;w)+\int_0^t\int_{\Rd} \delta^{\mathfrak{K}_z}(t-s,x,z;w)q(s,z,y)\,dzds\,.
\end{align*}
We deduce
\eqref{ineq:some-est_p_kappa-crit1} from 
\eqref{e:fract-der-est1-crit},
Lemma~\ref{ineq:I_0_oszagorne-crit1},
\cite[(92), (93)]{GS-2018}.
The inequality \eqref{ineq:some-est_p_kappa_1-crit1} results from
\eqref{ineq:some-est_gen_phi_xy-crit1},
\cite[Lemma~5.15]{GS-2018},
Lemma~\ref{l:some-estimates-3b-crit1}, \ref{l:convolution}(a)
and~\ref{lem:cal_TCh}.
\qed

\begin{lemma}\label{e:fract-der-p-kappa-2b-crit1}
Assume $\Qzero$ and $1-\lah <\beta\land \lah$.
Let $\beta_1\in (0,\beta]\cap (0,\lah)$.
 For every $T>0$ there exists a constant 
$c=c(d,T,\param,\kappa_2,\kb,\beta_1)$
such that for all $t\in (0,T]$, $x\in\Rd$,
\begin{equation*}
\left| \int_{\Rd}\nabla_x p^{\kappa}(t,x,y)\,dy\right|\leq  c \left[h^{-1}(1/t)\right]^{-1+\beta_1} \,, 
\end{equation*}
\end{lemma}
\pf
We get the inequality from Lemma~\ref{e:some-estimates-2bb-crit1}, \eqref{e:gradient-phi-y-crit1},
\eqref{eq:grad_phi_pomoc-crit1} and
 Lemma~\ref{lem:int_grad_phi}.
\qed

\begin{lemma}\label{l:p-kappa-difference-crit-1} 
Assume $\Qa$ or $\Qb$.\\
\noindent 
(a)
The function $p^{\kappa}(t,x,y)$ 
is jointly continuous   
on $(0, \infty)\times \Rd \times \Rd$.

\noindent
(b)  For every $T> 0$ there is a constant 
$c=c(d,T,\param,\kappa_2,\kb, \beta)$
 such that for all $t\in (0,T]$ and $x,y\in \Rd$,
$$
|p^{\kappa}(t,x,y)|\leq c t \err{0}{0}(t,x-y).
$$

\noindent
(c)  For all $t>0$, $x,y\in\Rd$, $x\neq y$,  
$$
\partial_t p^{\kappa}(t,x,y)= \LL_x^{\kappa}\, p^{\kappa}(t,x,y)\,.
$$

\noindent
(d)
For every $T>0$ there is a constant 
$c=c(d,T,\param,\kappa_2,\kb,\beta)$
 such that for all
$t\in (0,T]$, $x,y\in\Rd$,
\begin{align}\label{e:fract-der-p-kappa-1b-crit1}
|\LL_x^{\kappa} p^{\kappa}(t, x, y)|\leq c \err{0}{0}(t,x-y)\,,
\end{align}
and 
 \begin{align}\label{e:fract-der-p-kappa-2-crit1}
\left|\nabla_x p^{\kappa}(t,x,y)\right|\leq  c\! \left[h^{-1}(1/t)\right]^{-1} t \err{0}{0}(t,x-y)\,. 
\end{align}

\noindent
(e) For all $T>0$, $\gamma \in [0,1]\cap [0,\lah)$,
there is a constant  $c=c(d,T,\param,\kappa_2,\kb, \beta,\gamma)$
 such that for all $t\in (0,T]$ and $x,x',y\in \Rd$,
\begin{align*}
\left|p^{\kappa}(t,x,y)-p^{\kappa}(t,x',y)\right| \leq c 
 (|x-x'|^{\gamma}\land 1) \,t \left( \err{-\gamma}{0} (t,x-y)+ \err{-\gamma}{0}(t,x'-y) \right).
\end{align*}

\noindent
For all $T>0$, $\gamma \in [0,\beta)\cap [0,\lah)$,
there is a constant
$c=c(d,T,\param,\kappa_2,\kb, \beta,\gamma)$
such that for all $t\in (0,T]$ and $x,y,y'\in \Rd$,
\begin{align*}
\left|p^{\kappa}(t,x,y)-p^{\kappa}(t,x,y')\right| \leq c 
(|y-y'|^{\gamma}\land 1)\, t \left(  \err{-\gamma}{0}(t,x-y)+  \err{-\gamma}{0}(t,x-y') \right).
\end{align*}

\noindent
(f) The function $\LL_x^{\kappa}p^{\kappa}(t,x,y)$  is jointly continuous on $(0,\infty)\times \Rd\times \Rd$.
\end{lemma}
\pf 
The statement of (a) follows from Lemmas~\ref{lem:cont_frcoef} and~\ref{lem:phi_cont_joint-crit1}.
Part  (b) is a result of
Proposition~\ref{prop:gen_est_crit} 
and Lemma~\ref{lem:phi_cont_xy-crit1}.
 The equation  in (c) is a consequence of \eqref{e:p-kappa}, \eqref{eq:p_gen_klas} and
Corollary~\ref{e:phi-y-partial_1-crit}:
$\partial_t p^{\kappa}(t,x,y)=\LL_x^{\mathfrak{K}_x} p^{\mathfrak{K}_y}(t,x,y)+ \LL_x^{\mathfrak{K}_x} \phi_y(t,x)=\LL_x^{\mathfrak{K}_x} p^{\kappa}(t,x,y)$.
We get
\eqref{e:fract-der-p-kappa-1b-crit1} 
by 
\eqref{e:p-kappa},
\eqref{ineq:Lkp_abs},
\eqref{ineq:I_0_oszagorne-crit1-impr},
\cite[(92), (93)]{GS-2018}
(see also Lemma~\ref{e:L-on-phi-y-crit1} and \eqref{e:L-on-phi-y2-crit1}).
For the proof of \eqref{e:fract-der-p-kappa-2-crit1} we use Proposition~\ref{prop:gen_est_crit} and 
\eqref{e:gradient-phi-y-estimate-crit1}.
The first inequality of part  (e) follows from
 Lemmas~\ref{lem:pkw_holder} and~\ref{lem:phi_cont_xy-crit1}, and \cite[(92), (93)]{GS-2018}.
The same argument suffices for the second inequality of part (e) when supported by
$$
|p^{\mathfrak{K}_y}(t,x,y)-p^{\mathfrak{K}_{y'}}(t,x,y')|
\leq |p^{\mathfrak{K}_y}(t,-y,-x)-p^{\mathfrak{K}_{y}}(t,-y',-x)|
+ |p^{\mathfrak{K}_y}(t,x,y')-p^{\mathfrak{K}_{y'}}(t,x,y')|
$$
and 
Proposition~\ref{prop:Hcont_kappa_crit1}.
Part (f)
follows from 
Lemmas~\ref{lem:cont_frcoef} and~\ref{lem:Lphi_cont-crit}.
\qed

\section{Main Results and Proofs}\label{sec:Main}
In the whole section, we assume that either $\Qa$ or $\Qb$ holds.
\subsection{A nonlocal maximum principle}
Recall that $\LL^{\kappa,0^+}f:=\lim_{\varepsilon \to 0^+}\LL^{\kappa,\varepsilon}f$ is an extension of $\LL^{\kappa}f:=\LL^{\kappa,0}f$.
Moreover, 
the well-posedness of those operators requires the existence of the gradient $\nabla f$. 
The uniqueness of solutions to \eqref{e:nonlocal-max-principle-4-crit} stated in Corollary~\ref{cor:jedn_max-crit} will be used, for instance, in the next subsection.
For the proofs of the following, see \cite[Theorem~4.1]{GS-2018}.
\begin{theorem}
\label{t:nonlocal-max-principle-crit}
Let $T>0$ and   $u\in C([0,T]\times \Rd)$ be such that
\begin{align}\label{e:nonlocal-max-principle-1-crit}
\| u(t,\cdot)-u(0,\cdot) \|_{\infty} \xrightarrow {t\to 0^+} 0\,, \qquad \qquad \sup_{t\in [0,T]} \|  u(t,\cdot)\ind_{|\cdot|\geq r} \|_{\infty} \xrightarrow {r\to \infty}0\,.
\end{align}
Assume that  $u(t,x)$ satisfies the following equation: for all $(t,x)\in (0,T]\times \Rd$,
\begin{align}\label{e:nonlocal-max-principle-4-crit}
\partial_t u(t,x)=\LL_x^{\kappa,0^+}u(t,x)\, .
\end{align}
If $\sup_{x\in\Rd} u(0,x)\geq 0$, then
for every $t\in (0,T]$,
\begin{align}\label{e:nonlocal-max-principle-5-crit}
\sup_{x\in \R^d}u(t,x)\leq \sup_{x\in \Rd}u(0,x)\, .
\end{align}
\end{theorem}
\begin{corollary}\label{cor:jedn_max-crit}
If $u_1, u_2 \in C([0,T]\times \Rd)$ satisfy  \eqref{e:nonlocal-max-principle-1-crit}, \eqref{e:nonlocal-max-principle-4-crit} 
 and $u_1(0,x)=u_2(0,x)$, then $u_1\equiv u_2$ on $[0,T]\times \Rd$.
\end{corollary}

\subsection{Properties of the semigroup  $(P^{\kappa}_t)_{t\ge 0}$}
Define
$$
P_t^{\kappa}f(x)=\int_{\R^d}p^\kappa(t,x, y)f(y)dy.
$$
We first collect some properties of $\rr_t*f$.
\begin{remark}\label{rem:conv_Lp-crit}
We have $\rr_t*f \in C_b(\Rd)$
for any $f\in L^p(\Rd)$, $p\in [1,\infty]$.
Moreover,
$\rr_t*f\in C_0(\Rd)$ 
for any $f\in L^p(\Rd)\cup C_0(\Rd)$, $p\in [1,\infty)$.
Furthermore, there is $c=c(d)$ such that
$\|\rr_t*f \|_p\leq c \|f\|_p$ for all $t>0$, $p\in [1,\infty]$.
The above follows from $\rr_t\in L^1(\Rd)\cap L^{\infty}(\Rd)\subseteq L^q(\Rd)$ for every $q\in [1,\infty]$
(see \cite[Lemma~5.6]{GS-2018}),
and from properties of the convolution.
\end{remark}

\begin{lemma}\label{lem:bdd_cont-crit1}
(a)
We have $P_t^{\kappa} f \in C_b(\Rd)$ for any $f\in L^p(\Rd)$, $p\in [1,\infty]$.
Moreover, $P_t^{\kappa} f \in C_0(\Rd)$
for any $f\in L^p(\Rd)\cup C_0(\Rd)$, $p\in [1,\infty)$.
For every $T>0$ there exists a constant 
$c=c(d,T,\param,\kappa_2,\kb, \beta)$
such that for all $t\in(0,T]$ we get
$$
\|P^{\kappa}_t f\|_p\leq c \|f\|_p\,.
$$
(b) $P^{\kappa}_t\colon C_0(\Rd)\to C_0(\Rd)$, $t>0$, and for any bounded uniformly continuous function $f$,
$$
\lim_{t\to 0^+} \|P^{\kappa}_t f -f \|_{\infty}=0\,.
$$
(c)
$P^{\kappa}_t\colon L^p(\Rd)\to L^p(\Rd)$, $t>0$, $p\in [1,\infty)$, and for any $f\in L^p(\Rd)$,
$$
\lim_{t\to 0^+} \|P_t^{\kappa}f -f \|_p=0\,.
$$
\end{lemma}
\pf
The proof is like that for  \cite[Lemma~4.4]{GS-2018} and uses
Remark~\ref{rem:conv_Lp-crit}, 
Lemma~\ref{l:p-kappa-difference-crit-1},
\ref{lem:phi_cont_xy-crit1},
\ref{l:convolution}(a),
\eqref{e:some-estimates-2c-crit1},
\cite[Lemma~5.6]{GS-2018}
and
Proposition~\ref{prop:gen_est_crit}.
\qed

Our aim now is to prove
Proposition~\ref{lem:gen_sem_step1-crit1}
and
Lemma~\ref{lem:p-kappa-final-prop-crit1}.
The first one
provides an important link between $P_t^{\kappa}$ and the operator $\LL^{\kappa}$.
The other complements the fundamental properties of $p^\kappa(t,x, y)$.
They are both obtained by virtue of
Corollary~\ref{cor:jedn_max-crit}, but 
beforehand we have to make sure that certain functions satisfy
\eqref{e:nonlocal-max-principle-4-crit}.
The necessary results are prepared in Lemma~\ref{lem:grad_Pt-crit1} -- \ref{l:L-int-commute-crit1}.
In particular, we show that we can
apply 
the operator $\LL^{\kappa}$ to
$ \int_0^t P^{\kappa}_s f(x)\,ds$.

The formula \eqref{eq:grad_Pt_1-crit1} below is the same as \cite[(69)]{GS-2018}, however the proof given there fully relies on the condition $\lah>1$, while we assume that $\lah\leq 1$. Here we obtain \eqref{eq:grad_Pt_1-crit1} under an additional restriction on $f$, which suffices for our purposes.

\begin{lemma}\label{lem:grad_Pt-crit1}
For any $f\in L^p(\Rd)$, $p\in [1,\infty]$,  we have for all $t>0$, $x\in\Rd$,
\begin{align}\label{eq:grad_Pt-crit1}
\nabla_x \,P_t^{\kappa} f(x)= \int_{\Rd} \nabla_x\, p^{\kappa}(t,x,y) f(y)dy\,.
\end{align}
For any 
bounded (uniformly) H\"older continuous function
$f \in C^\eta_b(\Rd)$, $1-\lah<\eta$,
and all $t>0$, $x\in\Rd$,
\begin{align}\label{eq:grad_Pt_1-crit1}
\nabla_x \left( \int_0^t P^{\kappa}_s f(x)\,ds \right)= \int_0^t \nabla_x  P^{\kappa}_s f(x)\,ds\,.
\end{align}
\end{lemma}
\pf
 By  \eqref{e:fract-der-p-kappa-2-crit1} 
and 
\cite[Corollary~5.10]{GS-2018}
for $|\varepsilon|<h^{-1}(1/t)$,
\begin{align*}
 \left|  \frac1{\varepsilon}( p^{\kappa}(t,x+\varepsilon e_i,y)-p^{\kappa}(t,x,y)) \right| |f(y)| \leq c \left[h^{-1}(1/t)\right]^{-1}  \rr_t (x-y) |f(y)|\,.
\end{align*}
The right hand side is integrable by Remark~\ref{rem:conv_Lp-crit}. We can use the dominated convergence theorem, which gives
\eqref{eq:grad_Pt-crit1}.
 For $f \in C^\eta_b(\Rd)$ (we can assume that $\eta<\lah$)
 we let $\widetilde{x}=x+\varepsilon\theta e_i$ and by
\eqref{e:fract-der-p-kappa-2-crit1},
Lemma~\ref{e:fract-der-p-kappa-2b-crit1} and~\ref{l:convolution}(a) we have
\begin{align*}
&\left| \int_{\Rd}
  \frac1{\varepsilon}( p^{\kappa}(s,x+\varepsilon e_i,y)-p^{\kappa}(s,x,y)) f(y)\,dy \right|
\leq\left|  \int_{\Rd} \int_0^1 \partial_{x_i} p^{\kappa}(s,\widetilde{x},y) \, d\theta\, f(y)\,dy\right|\\
& \leq   \left| \int_{\Rd} \int_0^1  \partial_{x_i} p^{\kappa}(s,\widetilde{x},y)   \big[ f(y)-f(\widetilde{x})\big] \, d\theta \,dy\right| +  \left| \int_{\Rd} \int_0^1  \partial_{x_i} p^{\kappa}(s,\widetilde{x},y)f(\widetilde{x})\, d\theta \,dy\right|\\
&\leq c  \left[h^{-1}(1/s)\right]^{-1} \int_0^1 \int_{\Rd} s\err{0}{\eta}(s, \widetilde{x}-y)   \,dy\, d\theta 
+ c\left[h^{-1}(1/s)\right]^{-1+\beta_1}\\
&\leq c \left[h^{-1}(1/s)\right]^{-1+\eta}+  c \left[h^{-1}(1/s)\right]^{-1+\beta_1}\,.
\end{align*}
The right hand side is integrable over $(0,t)$ by 
\cite[Lemma~5.15]{GS-2018}.
 Finally, \eqref{eq:grad_Pt_1-crit1} follows by the dominated convergence theorem.
\qed

\begin{lemma}\label{l:L-int-commute0-crit1}
For any function $f\in L^p(\Rd)$, $p\in [1,\infty]$, and all $t>0$, $x\in\Rd$,
\begin{align}\label{e:L-int-commute-2-crit1}
\LL_x^{\kappa}P_t^{\kappa} f(x)=\int_{\Rd}\LL_x^{\kappa} \,p^{\kappa}(t,x, y)f(y)dy\, .
\end{align}
 Furthermore, 
for every $T>0$ there exists a constant $c>0$ such that
for all $f\in L^p(\Rd)$, $t\in (0,T]$, 
\begin{align}\label{e:LP-p-estimate-crit1}
\| \LL^{\kappa}P_t^{\kappa} f\|_p\leq c  t^{-1} \|f\|_p\,.
\end{align}
\end{lemma}
\pf
By the definition and \eqref{eq:grad_Pt-crit1},
\begin{align}\label{eq:LPf-crit1}
\LL_x^{\kappa} P_t^{\kappa} f(x)
=\int_{\Rd} \left(  \int_{\Rd} \delta^{\kappa}(t,x,y;z)  f(y)dy \right)  \kappa(x,z)J(z)dz\,.
\end{align}
The equality follows from an application of Fubini's theorem, justified by 
\eqref{ineq:some-est_p_kappa-crit1} and Remark~\ref{rem:conv_Lp-crit}.
The inequality follows then from \eqref{e:L-int-commute-2-crit1},  \eqref{e:fract-der-p-kappa-1b-crit1},
Remark~\ref{rem:conv_Lp-crit}.
\qed

\begin{lemma}\label{lem:for_max-crit1}
Let $f\in C_0(\Rd)$. For $t>0$, $x\in\Rd$ we define
$u(t,x)=P^{\kappa}_t f(x)$ and $u(0,x)=f(x)$. 
Then $u\in C([0,T]\times \Rd)$,
\eqref{e:nonlocal-max-principle-1-crit} holds and $\partial_t u(t,x)=\LL_x^{\kappa}u(t,x)$ for all $t,T>0$, $x\in\Rd$.
\end{lemma}
\pf
The proof is exactly like that of \cite[Lemma~4.7]{GS-2018}.
\qed

\begin{lemma}\label{l:L-int-commute-crit1}
For any bounded (uniformly) H\"older continuous function 
$f \in C^\eta_b(\Rd)$, $1-\lah<\eta$,
and all $t>0$, $x\in\Rd$,
we have $\int_0^t | \LL_x^{\kappa} P_s^{\kappa}f(x)|ds <\infty$ and
\begin{align}\label{e:L-int-commute-crit1}
\LL_x^{\kappa}\left( \int_0^t P_s^{\kappa}f(x)\,ds\right) =\int_0^t \LL_x^{\kappa} P_s^{\kappa}f(x)\,ds\,.
\end{align}
\end{lemma}
\pf
By definition and Lemma~\ref{lem:grad_Pt-crit1},
\begin{align*}
\LL_x^{\kappa} \int_0^t P_s^{\kappa}f(x)\,ds
&=\int_{\Rd} \left( \int_0^t \int_{\Rd} \delta^{\kappa} (s,x,y;z)   f(y)dy ds \right) \kappa(x,z)J(z)dz\,.
\end{align*}
Note that the proof will be finished by \eqref{eq:LPf-crit1} if we can change the order of integration from $dsdz$ to $dzds$. To this end we use Fubini's theorem, which is justified by the following. We have $|f(y)-f(x)|\leq c (|y-x|^{\eta} \land 1)$ and we can assume that $\eta<\lah$.
Then
\begin{align*}
\int_{\Rd}  \int_0^t &\left| \int_{\Rd} \delta^{\kappa} (s,x,y;z)   f(y)dy \right| ds \, \kappa(x,z)J(z)dz\\
&\leq \int_{\Rd}  \int_0^t \left| \int_{\Rd}  \delta^{\kappa} (s,x,y;z) \big[f(y)-f(x)\big] dy\right| ds \,\kappa(x,z)J(z)dz\\
&\quad+\int_{\Rd}  \int_0^t \left| \int_{\Rd} \delta^{\kappa} (s,x,y;z) f(x)  dy\right| ds\, \kappa(x,z)J(z)dz=: {\rm I}_1+{\rm I}_2\,.
\end{align*}
By \eqref{ineq:some-est_p_kappa-crit1}, we get
${\rm I}_1\leq c \int_0^t \int_{\Rd} \ct(s) \err{0}{\eta}(s,y-x)  dyds$, while by \eqref{ineq:some-est_p_kappa_1-crit1}
${\rm I}_2\leq c \int_0^t \ct(s) s^{-1} \left[h^{-1}(1/s)\right]^{\beta_1}ds$.
The integrals are finite by 
\cite[Lemma~5.15]{GS-2018},
Lemma~\ref{l:convolution}(a)
and~\ref{lem:cal_TCh}.
\qed

\begin{proposition}\label{lem:gen_sem_step1-crit1}
For any $f\in C_b^{2}(\Rd)$ and all $t>0$, $x\in\Rd$,
\begin{align}\label{eq:gen_sem_step1-crit1}
P_t^{\kappa}f(x)-f(x)=\int_0^t P_s^{\kappa}\LL^{\kappa} f(x)\,ds\,.
\end{align}
\end{proposition}
\pf
We outline the main steps, for details, see the proof of \cite[Proposition~4.9]{GS-2018}.

\noindent
(i) Note that $\LL^{\kappa}f \in C_0(\Rd)$ for any $f\in C_0^2(\Rd)$.

\noindent
(ii) Show that 
$\LL^{\kappa}f  \in C_b^{\eta}(\Rd)$
for any $f\in C_0^{2,\eta}(\Rd)$.

\noindent
(iii) 
Show that \eqref{eq:gen_sem_step1-crit1} holds for any $f\in C_0^{2,\eta}(\Rd)$ if $1-\lah<\eta\leq \beta$.
It is achieved by
using Corollary~\ref{cor:jedn_max-crit}
to prove that the following functions are equal
\begin{align*}
u_1(t,x)=
\begin{cases}
P_t^{\kappa} f(x), \quad &t>0\\
f(x), & t=0
\end{cases}
\,,
\qquad
u_2(t,x)=
\begin{cases}
f(x)+\int_0^t P_s^{\kappa}\LL^{\kappa} f(x)\,ds, \quad &t>0 \\
f(x), &t=0
\end{cases}
\,.
\end{align*}
Here we use Lemmas~\ref{lem:for_max-crit1}, \ref{l:L-int-commute-crit1} and \cite[Theorem~7.21]{MR924157}.

\noindent
(iv) We extend \eqref{eq:gen_sem_step1-crit1} to  $f\in C_b^2(\Rd)$ by approximating it with $f_n=(f*\phi_n)\cdot \varphi_n \in C_c^{\infty}(\Rd)$, where $\phi_n$ is a standard mollifier while $\varphi_n(x)=\varphi(x/n)$ for $\varphi\in C_c^{\infty}(\Rd)$ satisfying $\varphi(x)=1$ if $|x|\leq 1$, and $\varphi(x)=0$ if $|x|\geq 2$.
\qed

\begin{lemma}\label{lem:p-kappa-final-prop-crit1}
The function $p^{\kappa}(t,x,y)$ is non-negative, $\int_{\Rd} p^{\kappa}(t,x,y)dy= 1$ and $p^{\kappa}(t+s,x,y)=\int_{\Rd}p^{\kappa}(t,x,z)p^{\kappa}(s,z,y)dz$
for all $s,t>0$, $x,y\in\Rd$.
\end{lemma}
\pf
The proof is like that of \cite[Lemma~4.10]{GS-2018} and we use Lemma~\ref{lem:for_max-crit1}, Theorem~\ref{t:nonlocal-max-principle-crit}, Lemma~\ref{l:p-kappa-difference-crit-1}, Corollary~\ref{cor:jedn_max-crit}, Proposition~\ref{lem:gen_sem_step1-crit1}.
However, the proof of the convolution property in \cite{GS-2018} contains a gap: 
at that stage it is not clear why the function $p^{\kappa}(t+s,x,y)$ should satisfy the equation \eqref{e:nonlocal-max-principle-4-crit} for all $x\in\Rd$. Here we present a  correction that is valid for both papers.

Let $T,s>0$ and $\varphi\in C_c^{\infty}(\Rd)$.
For $t>0$, $x\in\Rd$ define
$$
u_1(t,x)=P_t^{\kappa}f(x)\,,\qquad u_1(0,x)=f(x)=P_s^{\kappa}\varphi(x)\,,
$$
and
$$
u_2(t,x)=P_{t+s}^{\kappa}\varphi(x)\,, \qquad u_2(0,x)=P_s^{\kappa}\varphi(x)\,.
$$
By Lemma~\ref{lem:bdd_cont-crit1}(b) $f\in C_0(\Rd)$
and thus by Lemma~\ref{lem:for_max-crit1}
$u_1$ satisfies the assumptions of Corollary~\ref{cor:jedn_max-crit}. Now, since $\varphi$
has compact support by Lemma~\ref{l:p-kappa-difference-crit-1}(a) we get $u_2\in C([0,T]\times \Rd)$. We will use 
\cite[(94)]{GS-2018}
several times in what follows. By Lemma~\ref{l:p-kappa-difference-crit-1} (c) and (d),
\begin{align*}
\|u_2(t,\cdot)-u_2(0,\cdot)\|_{\infty}
&\leq \sup_{x\in\Rd}\int_{\Rd}\int_0^t |\LL_x^{\kappa}\,p^{\kappa}(u+s,x,y)|du \,|\varphi(y)|dy\\
&\leq c t  \err{0}{0}(s,0)\int_{\Rd}|\varphi(y)|dy
\to 0\,, \quad \mbox{as } t\to 0^+\,.
\end{align*}
Furthermore, by Lemma~\ref{l:p-kappa-difference-crit-1}(b)
\begin{align*}
\sup_{t\in[0,T]}\|u_2(t\cdot)\ind_{|\cdot|\geq r}\|_{\infty}
&\leq c (T+s) \sup_{x\in\Rd} \ind_{|x|\geq r}\int_{\Rd}\err{0}{0}(s,x-y)|\varphi(y)|dy\\
&=c (T/s+1)\sup_{x\in\Rd}  \ind_{|x|\geq r} \left(\rr_t*|\varphi|\right)(x)\to 0\,,
\quad \mbox{as } r\to \infty\,,
\end{align*}
because 
$\rr_t*|\varphi|\in C_0(\Rd)$ (see Remark~\ref{rem:conv_Lp-crit}).
Finally, 
by the mean value theorem, Lemma~\ref{l:p-kappa-difference-crit-1}(c), \eqref{e:fract-der-p-kappa-1b-crit1} and the dominated convergence theorem $\partial_t u_2(t,x)=\int_{\Rd}\partial_t p^{\kappa}(t+s,x,y)\varphi(y) dy$. Then we apply Lemma~\ref{l:p-kappa-difference-crit-1}(c) and Lemma~\ref{l:L-int-commute0-crit1} to obtain $\partial_t u_2(t,x)=\LL_x^{\kappa}u(t,x)$. Therefore, by Corollary~\ref{cor:jedn_max-crit} $u_1=P_t^{\kappa}P_s^{\kappa}\varphi = P_{t+s}^{\kappa}\varphi = u_2$. The convolution property now follows by Fubini theorem, 
the arbitrariness of $\varphi$ and by Lemma~\ref{l:p-kappa-difference-crit-1}(a).
\qed

\subsection{Proofs of Theorems~\ref{t:intro-main}--\ref{thm:onC0Lp}}

At this point we have all necessary tools to proceed exactly like in \cite{GS-2018}. A slight difference in the proof of Theorem~\ref{thm:lower-bound} is explained below.\\

\noindent
{\bf Proof of Theorem~\ref{thm:onC0Lp}}.
It is the same as in \cite{GS-2018} and relies on
Lemma~\ref{lem:bdd_cont-crit1}, \ref{lem:p-kappa-final-prop-crit1},  Proposition~\ref{lem:gen_sem_step1-crit1},
\eqref{eq:grad_Pt-crit1}, \eqref{e:fract-der-p-kappa-2-crit1}, Remark~\ref{rem:conv_Lp-crit},
\eqref{ineq:some-est_p_kappa-crit1}, \eqref{eq:LPf-crit1},
\eqref{e:intro-kappa},
\eqref{e:intro-kappa-holder},
\eqref{e:L-int-commute-2-crit1},  Lemma~\ref{l:p-kappa-difference-crit-1}(f) and~(d),
\cite[Corollary~5.10]{GS-2018},
 \eqref{e:LP-p-estimate-crit1},
\cite[Chapter~1, Theorem~2.4(c) and~2.2]{MR710486},
\cite[Chapter~2, Theorem~5.2(d)]{MR710486}.
\qed

\noindent
{\bf Proof of Theorem~\ref{t:intro-further-properties}}.
All the results are collected in
Lemma~\ref{l:p-kappa-difference-crit-1}
and~\ref{lem:p-kappa-final-prop-crit1}, except for part (8), which is given in Theorem~\ref{thm:onC0Lp} part 3(c).
\qed

\noindent
{\bf Proof of Theorem~\ref{t:intro-main}}.
The same as in \cite{GS-2018}.
\qed

\noindent
{\bf Proof of Theorem~\ref{thm:lower-bound}}.
The argument is the same as in \cite{GS-2018} except for the following modification of the prove of the estimate
$$
\sup_{|\xi|\leq 1/r} |q(z,\xi)|
\leq  c_2 h(r).
$$
Since for
$\varphi \in \RR$ we have $\left|e^{i\varphi}-1-i\varphi \right|\leq  |\varphi|^2$,
by \eqref{e:intro-kappa-crit} we get
\begin{align*}
|q(z,\xi)|&\leq
\int_{\Rd} \left| e^{i\left<\xi, w\right>}-1 - i\left<\xi,w\right> \ind_{|w|<1\land 1/|\xi|} \right|\kappa(x,w)J(w)dw\\
&\quad +\left| i \left<\xi,\int_{\Rd} w \left(\ind_{|w|<1\land 1/|\xi|}
- \ind_{|w|<1}\right)
\kappa(x,w)J(w)dw \right>\right|\\
&\leq
|\xi|^2 \int_{|w|<1\land 1/|\xi|} |w|^2 \kappa(x,w)J(w)dw
+ \int_{|w|\geq 1\land 1/|\xi|} 2\, \kappa(x,w)J(w)dw\\
&\quad + |\xi| \left| \int_{\Rd} w \left(\ind_{|w|<1\land 1/|\xi|}
- \ind_{|w|<1}\right)
\kappa(x,w)J(w)dw\right|
\leq c h(1\land 1/|\xi|)\,.
\end{align*}
\qed

\section{Appendix - Unimodal L{\'e}vy processes}\label{sec:appA}

Let $d\in\N$ and
$\nu:[0,\infty)\to[0,\infty]$ be a non-increasing  function  satisfying
$$\int_{\Rd}  (1\land |x|^2) \nu(|x|)dx<\infty\,.$$
For any such $\nu$, there exists a unique
 pure-jump
isotropic unimodal L{\'e}vy process $X$ (see \cite{MR3165234}, \cite{MR705619}).
We define
$h(r)$, $K(r)$ and $\rr_t(x)$ as in the introduction.
At this point we refer the reader to \cite[Section~5]{GS-2018}
for various important properties of 
those functions. 
Following \cite[Section~5]{GS-2018}
in the whole section {\bf we assume that} $h(0^+)=\infty$.
We consider the scaling conditions: 

\noindent
there are $\lah \in(0,2]$, $C_h\in[1,\infty)$ and $\theta_h\in(0,\infty]$ such that
\begin{equation}\label{eq:wlsc:h}
 h(r)\leq C_h\lambda^{\lah }h(\lambda r),\qquad \lambda\leq 1,\, r< \theta_h;
 \end{equation} 
there are $\uah \in (0,2]$, $c_h\in (0,1]$ and $\theta_h \in (0,\infty]$ 
such that
\begin{equation}\label{eq:wusc:h}
 c_h\,\lambda^{\uah}\,h(\lambda r)\leq h(r)\, ,\quad \lambda\leq 1, \,r< \theta_h.
\end{equation}
 
The first and the latter inequality in the next lemma are taken from \cite[Section~5]{GS-2018}. We keep them here for easier reference.

\begin{lemma}\label{lem:int_J}
Let $h$ satisfy \eqref{eq:wlsc:h} with $\lah>1$, then
\begin{align*}
\int_{r \leq |z|<  \theta_h }  |z|\nu(|z|)dz \leq \frac{(d+2) C_h}{\lah-1} \, r h(r)\,, \qquad r>0\,.
\end{align*}
Let $h$ satisfy \eqref{eq:wlsc:h} with $\lah=1$, then
\begin{align*}
\int_{r \leq |z|<  \theta_h }  |z|\nu(|z|)dz \leq [(d+2)C_h] \, \ln(\theta_h/r)\, r h(r)\,, \qquad r>0\,.
\end{align*}
Let $h$ satisfy \eqref{eq:wusc:h} with $\uah<1$, then
\begin{align*}
\int_{|z|< r} |z| \nu(|z|)dz\leq \frac{d+2}{c_h(1-\uah)}\, r h(r)\,,\qquad r< \theta_h \,.
\end{align*}
\end{lemma}
\pf
Under \eqref{eq:wlsc:h} with $\lah=1$ we have
\begin{align*}
\int_{r \leq |z|<  \theta_h }  |z|\nu(|z|)dz 
\leq (d+2) \int_r^{\theta_h} h(s)ds 
\leq (d+2) C_h\int_r^{\theta_h} (r/s)h(r)ds\,,
\end{align*}
which ends the proof.
\qed

\begin{lemma}\label{lem:cal_TCh}
Assume  \eqref{eq:intro:wlsc}.
Let $k,l \geq 0$ and $\theta, \eta, \beta,\gamma\in\R$ satisfy $(\beta/2)\land (\beta/\lah)+1-\theta>0$, $(\gamma/2)\land(\gamma/\lah)+1-\eta>0$. 
For every $T>0$ 
there  exists a constant $c=c(\lah,C_h,h^{-1}(1/T)\vee 1, \theta, \eta,\beta,\gamma,k,l)$ such that for all $t\in (0,T]$,
\begin{align*}
\int_0^t [\TCh(u)]^{l}\, u^{-\eta}\left[ h^{-1}(1/u)\right]^{\gamma}[\TCh(t-u)]^{k}\, (t-u)^{-\theta}\left[ h^{-1}(1/(t-u))\right]^{\beta}du \\
\leq c\, [\TCh(t)]^{l+k}\, t^{1-\eta-\theta}\left[h^{-1}(1/t)\right]^{\gamma+\beta}\,.
\end{align*}
Furthermore, $\TCh(t/2)\leq c\, \TCh(t)$, $t\in (0,T]$.
\end{lemma}
\pf
The last part of the statement follows from  
\cite[Lemma~5.3 and Remark~5.2]{GS-2018}.
Note that it suffices to consider the integral over $(0,t/2)$.
Again by
\cite[Lemma~5.3 and Remark~5.2]{GS-2018}
we have for
$c_0=C_h [h^{-1}(1/T)\vee 1]^2$ and $s\in (0,1)$,
$$
\left[ h^{-1}(t^{-1}s^{-1})\right]^{-1}
\leq c_0\, s^{-1/\lah}\left[h^{-1}(1/t)\right]^{-1}, \qquad
h^{-1}(t^{-1}s^{-1})\leq  s^{1/2}h^{-1}(t^{-1})\,.
$$
Thus for $u\in (0,t/2)$ we get
$$ 
[\TCh(t-u)]^{k}(t-u)^{-\theta}[h^{-1}(1/(t-u))]^{\beta}\leq c \,[\TCh(t)]^{k}\, t^{-\theta}\,[h^{-1}(1/t)]^{\beta}\,,
$$
and we concentrate on 
$$\int_0^{t/2} [\TCh(u)]^{l} u^{-\eta} \left[h^{-1}(1/u)\right]^{\gamma}du
\leq 
c\,t^{1-\eta} \left[h^{-1}(1/t)\right]^{\gamma} \int_0^{1/2}  [\TCh(ts)]^{l} s^{(\gamma/2)\land (\gamma/\lah)-\eta} ds
\,.$$
Furthermore, we have
\begin{align*}
&\int_0^{1/2} \left[ \ln\left(1\vee\left[ h^{-1}(t^{-1}s^{-1})\right]^{-1}\right)\right]^{l} s^{(\gamma/2)\land (\gamma/\lah)-\eta}ds\\
&\leq \int_0^{1/2} 2^{l} \left\{ [\ln (c_0 s^{-1})]^{l} +\left[ \ln \left(1\vee\left[ h^{-1}(t^{-1}\right]^{-1}\right)\right]^{l} \right\} s^{(\gamma/2)\land (\gamma/\lah)-\eta}ds
\leq  c \,[\TCh(t)]^{l}\,.
\end{align*}
Finally,
\begin{align*}
 \int_0^{1/2}  [\TCh(ts)]^{l} s^{(\gamma/2)\land (\gamma/\lah)-\eta} ds
\leq c\, [\TCh(t)]^{l}.
\end{align*}
\qed

The next lemma is taken from 
\cite[Section~5]{GS-2018} and complemented with part (d).
It is one of the most often used technical
result in the paper.
Let $B(a,b)$ be the beta function, i.e., $B(a,b)=\int_0^1 s^{a-1} (1-s)^{b-1}ds$, $a,b>0$.

\begin{lemma}\label{l:convolution}
Assume \eqref{eq:intro:wlsc} and let $\beta_0\in(0,\lah\land 1)$.
\begin{itemize}
\item[(a)] For every $T>0$ there exists a constant $c_1=c_1(d,\beta_0,\lah,C_h,h^{-1}(1/T)\vee 1)$ such that for all $t\in(0,T]$ and $\beta\in [0,\beta_0]$,
$$
\int_{\Rd} \err{0}{\beta}(t,x)\,dx \leq c_1 t^{-1} \left[h^{-1}(1/t)\right]^{\beta}\,.
$$
\item[(b)] 
For every $T>0$ there exists a constant $c_2=c_2(d,\beta_0,\lah,C_h,h^{-1}(1/T)\vee 1) \ge 1$
such that
for all $\beta_1,\beta_2,n_1,n_2,m_1,m_2 \in[0,\beta_0]$ with $n_1, n_2 \leq \beta_1+\beta_2$, $m_1\leq \beta_1$, $m_2\leq \beta_2$ and all $0<s<t\leq T$, $x\in\Rd$,
\begin{align*}
\int_{\Rd} \err{0}{\beta_1}(t-s&,x-z)\err{0}{\beta_2}(s,z) \,dz\\
\leq  
c_2 &\Big[ \left( (t-s)^{-1} \left[h^{-1}(1/(t-s))\right]^{n_1} + s^{-1}\left[h^{-1}(1/s)\right]^{n_2}\right) \err{0}{0}(t,x)\\
&+(t-s)^{-1}\left[ h^{-1}(1/(t-s))\right]^{m_1} \err{0}{\beta_2}(t,x) + 
s^{-1}\left[ h^{-1}(1/s)\right]^{m_2} \err{0}{\beta_1}(t,x)\Big].
\end{align*}
\item[(c)] Let $T>0$. For all 
$\gamma_1, \gamma_2\in\RR$, 
$\beta_1,\beta_2,n_1,n_2,m_1,m_2 \in[0,\beta_0]$ with $n_1, n_2 \leq \beta_1+\beta_2$, $m_1\leq \beta_1$, $m_2\leq \beta_2$
and $\theta,\eta \in [0,1]$,
satisfying 
\begin{align*}
(\gamma_1+n_1\land m_1)/2 \land (\gamma_1+n_1\land m_1)/\lah +1-\theta>0\,,\\
(\gamma_2+n_2\land m_2)/2\land (\gamma_2+n_2\land m_2)/\lah+1-\eta>0\,,
\end{align*}
and all
$0<s<t\leq T$, $x\in\Rd$, we have
\begin{align}
\int_0^t\int_{\Rd}& (t-s)^{1-\theta}\, \err{\gamma_1}{\beta_1}(t-s,x-z) \,s^{1-\eta}\,\err{\gamma_2}{\beta_2}(s,z) \,dzds  \nonumber\\
& \leq c_3 \,
t^{2-\eta-\theta}\Big( \err{\gamma_1+\gamma_2+n_1}{0}
+\err{\gamma_1+\gamma_2+n_2}{0}
+\err{\gamma_1+\gamma_2+m_1}{\beta_2} 
+\err{\gamma_1+\gamma_2+m_2}{\beta_1} 
\Big)(t,x)\,,\label{e:convolution-3}
\end{align}
where 
$
c_3= c_2 \, (C_h[h^{-1}(1/T)\vee 1]^2)^{-(\gamma_1\land 0+\gamma_2 \land 0)/\lah}
 B\left( k+1-\theta, \,l+1-\eta\right)
$
and 
\begin{align*}
k=\left(\frac{\gamma_1+n_1\land m_1}{2}\right)\land
\left(\frac{\gamma_1+n_1\land m_1}{\lah}\right),
\quad l=\left(\frac{\gamma_2+n_2\land m_2}{2}\right)\land \left(\frac{\gamma_2+n_2\land m_2}{\lah}\right).
\end{align*}
\item[(d)] Let $T>0$.
For all $k,l\geq 0$,
$\gamma_1, \gamma_2\in\RR$, 
$\beta_1,\beta_2,n_1,n_2,m_1,m_2 \in[0,\beta_0]$ with $n_1, n_2 \leq \beta_1+\beta_2$, $m_1\leq \beta_1$, $m_2\leq \beta_2$
and $\theta,\eta \in [0,1]$,
satisfying 
\begin{align*}
(\gamma_1+n_1\land m_1)/2 \land (\gamma_1+n_1\land m_1)/\lah +1-\theta>0\,,\\
(\gamma_2+n_2\land m_2)/2\land (\gamma_2+n_2\land m_2)/\lah+1-\eta>0\,,
\end{align*}
and for all 
$0<s<t\leq T$, $x\in\Rd$, we have
\begin{align}
\int_0^t\int_{\Rd}& [\TCh(t-s)]^k (t-s)^{1-\theta}\, \err{\gamma_1}{\beta_1}(t-s,x-z) [\TCh(s)]^{l}\,s^{1-\eta}\,\err{\gamma_2}{\beta_2}(s,z) \,dzds  \nonumber\\
& \leq c_4 \, [\TCh(t)]^{k+l}\,
t^{2-\eta-\theta}\Big( \err{\gamma_1+\gamma_2+n_1}{0}
+\err{\gamma_1+\gamma_2+n_2}{0}
+\err{\gamma_1+\gamma_2+m_1}{\beta_2} 
+\err{\gamma_1+\gamma_2+m_2}{\beta_1} 
\Big)(t,x)\,,\label{e:convolution-3-crit}
\end{align}
where $c_4=c_4(d,\beta_0,\lah,C_h,h^{-1}(1/T)\vee 1,k,l,\gamma_1,\gamma_2,n_1,n_2,m_1,m_2,\theta,\eta)$.
\end{itemize}
\end{lemma}
\pf
For the proof of part (d) we multiply the result of part (b) by
$$
[\TCh(t-s)]^{k} (t-s)^{1-\theta} \left[h^{-1}(1/(t-s))\right]^{\gamma_1} [\TCh(s)]^{\gamma_1} s^{1-\eta} \left[h^{-1}(1/s)\right]^{\gamma_2}\,,
$$
and apply Lemma~\ref{lem:cal_TCh}.
\qed

\begin{remark}
When using Lemma~\ref{l:convolution} without specifying the parameters, we apply the usual case, i.e., $n_1=n_2=\beta_1+\beta_2$ ($\leq \beta_0$), $m_1=\beta_1$, $m_2=\beta_2$. Similarly, if only $n_1$, $n_2$ are specified, then $m_1=\beta_1$, $m_2=\beta_2$.
\end{remark}

\section*{Acknowledgements}

The author would like to thank Krzysztof Bogdan and Tomasz Jakubowski for discussions and comments that helped to improve the clarity of the manuscript.

\small
\bibliographystyle{abbrv}

\end{document}